\date{}

\documentclass[12pt,reqno]{amsart}
\usepackage{latexsym,amsmath,amsfonts,amscd,amssymb}
\usepackage{graphics}
\usepackage{float}

\theoremstyle{plain}
\newtheorem{theorem}{Theorem}[section]

\newtheorem{lemma}[theorem]{Lemma}

\theoremstyle{definition}
\newtheorem{definition}[theorem]{Definition}

\newtheorem{conjecture}[theorem]{Conjecture}
\numberwithin{equation}{section}

\newcommand{\g}{\gamma}

\newcommand{\n}{\nu}

\newcommand{\bk}{\mathbf{k}}

\newcommand{\cA}{\mathcal{A}}

\newcommand{\cK}{\mathcal{K}}

\newcommand{\cL}{\mathcal{L}}
\newcommand{\cO}{\mathcal{O}}
\newcommand{\cP}{\mathcal{P}}
\newcommand{\cR}{\mathcal{R}}
\newcommand{\cS}{\mathcal{S}}

\renewcommand{\AA}{\mathbb{A}}
\newcommand{\FF}{\mathbb{F}}
\newcommand{\QQ}{\mathbb{Q}}

\newcommand{\RR}{\mathbb{R}}
\newcommand{\CC}{\mathbb{C}}
\newcommand{\HH}{\mathbb{H}}

\newcommand{\NN}{\mathbb{N}}
\newcommand{\PP}{\mathbb{P}}
\newcommand{\TT}{\mathbb{T}}
\newcommand{\ZZ}{\mathbb{Z}}

\newcommand{\abs}[1]{\lvert#1\rvert}

\newcommand{\Gal}{\mathrm{Gal}}
\newcommand{\Li}{\mathrm{li}}
\newcommand{\Tr}{\mathrm{Tr}}
\newcommand{\Det}{\mathrm{Det}}

\title{Notes on the Riemann Hypothesis}


\author[R. P\'{e}rez-Marco]{Ricardo P\'{e}rez-Marco}
\address{CNRS, IMJ-PRG, Paris 7, Bo\^\i te courrier 7012, 75005 Paris Cedex 13, France}
\email{ricardo.perez.marco@gmail.com}

\thanks{We added section 3.2 and 3.5 (February 2018). These notes grew originally from a set of conferences given on March 2010 
on the ``Problems of the Millenium`` at 
the Complutense University in Madrid. I am grateful to the organizers Ignacio Sols and Vicente Mu\~noz.  
I also thank  Javier Fres\'an 
for his assistance with the notes from my talk, Jes\'us Mu\~noz, Arturo \'Alvarez, Santiago Boza and Javier Soria for their 
reading and comments on earlier versions of this text. I thank also Wemp Pacheco and Jan Heering for pointing out 
corrections in previous versions of the manuscript. Finally I also thank 
the Real Sociedad Matem\'atica 
Espa\~nola (R.S.M.E.) and the organizers of 
the event on the presentation of the Millenium problems celebrating the centenary anniversary of the RSME for 
the opportunity to participate in this event and present this subject.}

\begin{document}

\begin{abstract}
Our aim is to give an introduction to the Riemann Hypothesis and a panoramic view of the world of 
zeta and L-functions. 
We first review Riemann's foundational article and discuss the mathematical background of 
the time and his possible motivations for making his famous conjecture. We discuss some of the most relevant developments 
after Riemann that have contributed to a better understanding of the conjecture.
\end{abstract}

\maketitle

\tableofcontents

\section{Euler transalgebraic world.}

The first non-trivial occurrence of the real zeta function 
\begin{equation*}
\zeta(s)=\sum_{n=1}^\infty {n^{-s}}, \quad s \in \mathbb{R}, \ \ s>1
\end{equation*} 
appears to be in 1735 when Euler computed the infinite sum of the inverse of squares, the so called ``Basel Problem", a problem originally raised by P. Mengoli and  studied by Jacob Bernoulli,
\begin{equation*}
\zeta(2)=\sum_{n=1}^\infty {n^{-2}} \ .
\end{equation*}
Euler's solution was based on the {\it transalgebraic} interpretation of $\zeta(2)$ as the Newton sum of exponent $s$, for $s=2$, of 
the transcendental function $\frac{\sin\pi z}{\pi z}$, which vanishes at each non-zero integer. As Euler would write:
\begin{equation*}
\frac{\sin\pi z}{\pi z}=\prod_{n \in \ZZ^*} \left (1-\frac{z}{n} \right ) \ .
\end{equation*} 
Or, more precisely, regrouping the product for proper convergence (in such a way that respects an underlying symmetry, 
here $z\mapsto -z$, is a recurrent important idea),
\begin{eqnarray*}
\frac{\sin\pi z}{\pi z}&&=\left (1-z \right )\left (1+z\right )\left (1-\frac{z}{2}\right )\left (1+\frac{z}{2}\right )\ldots  
=\prod_{n=1}^\infty \left (1-\frac{z^2}{n^2}\right )\\
&&=1-\left (\sum_{n=1}^\infty \frac{1}{n^2}\right )z^2+\left (\sum_{n,m=1}^\infty \frac{1}{n^2 m^2}\right )z^4+ \ldots \\
\end{eqnarray*} 
On the other hand, the power series expansion of the sine function gives  
\begin{eqnarray*}
\frac{\sin\pi z}{\pi z}&&=\frac{1}{\pi z}\left (\pi z-\frac{1}{3!}(\pi z)^3+\frac{1}{5!}(\pi z)^5-\ldots \right ) \\
&&=1-\frac{\pi^2}{6}\ z^2+\frac{\pi^4}{120}\ z^4-\ldots,\\
\end{eqnarray*} 
and identifying coefficients, we get
\begin{equation*}
\zeta(2)=\frac{\pi^2}{6} \ .
\end{equation*}
In the same way, one obtains 
$$
\zeta(4)=\frac{\pi^4}{90} \ ,
$$
as well as the values $\zeta(2k)$ for each positive integer value of $k$. This seems to be 
one of the first occurrences of symmetric functions in an infinite number of variables, something that is at the heart of Transalgebraic Number Theory.

\medskip

This splendid argument depends crucially on the fact that the only zeros of $\sin(\pi z)$ in the 
complex plane are on the real line, and therefore are reduced to the natural integers. This, at the time, 
was the central point for a serious objection: It was necessary to prove the non-existence 
of non-real zeros of $\sin (\pi z)$. This follows easily from the unveiling of the relationship of trigonometric functions 
with the complex exponential (due to Euler) or to Euler $\Gamma$-function via the complement (or reflection) formula (due to Euler). 
Thus, we already meet at this early stage the importance of establishing that the divisor of a transcendental 
function lies on a line. We highlight this notion:

\medskip

 \textbf{Divisor on lines property (DL).} {\it A meromorphic function on the complex plane has its divisor on lines, 
or has the \textbf{DL} property, if its divisor is contained in a finite number of horizontal and vertical lines.}

\medskip

Meromorphic functions satisfying the \textbf{DL} property do form a multiplicative group. 
The \textbf{DL} property is close to having 
an order two symmetry or satisfying a functional equation. More precisely, a meromorphic function of finite order satisfying 
the \textbf{DL} property with its divisor contained in a single line, after multiplication by the exponential of a polynomial, 
does satisfy a functional equation that corresponds to the symmetry with respect to the divisor line. 
Obviously the converse is false. In our example the symmetry discussed comes from real analyticity.

\medskip

Another transcendental function, also well known to Euler, satisfying the \textbf{DL} property, is Euler gamma function,
\begin{equation*}
\Gamma(s)=\lim_{n\rightarrow \infty} \frac{n^s n!}{s(s+1)\cdots (s+n)} \ . 
\end{equation*}
It interpolates the factorial function and has its divisor, consisting only of poles, at the non-positive integers. 
Again, real analyticity gives one functional equation. The zeta function $\zeta(s)$ is exactly Newton sum of power $s$ 
for the zeros of the entire function $\Gamma^{-1}$, which is given by its Weierstrass product 
(where $\gamma=0.5772157\ldots$ is Euler constant)
$$
\frac{1}{\Gamma (s)}= s e^{\gamma s} \prod_{n=1}^{+\infty } \left (1+\frac{s}{n}\right ) \ e^{-s/n} \ .
$$
The $\Gamma$-function decomposes in the \textbf{DL} group our previous sine function by the complement formula
$$
\frac{\sin\pi z}{\pi z}=\frac{1}{\Gamma (1+z)} \ .\  \frac{1}{\Gamma (1-z)} \ .
$$

\bigskip 

This ``interpolation idea" which is at the origin of the definition of the $\Gamma$-function 
permeates the work of 
Euler and is also central in the ideas of Ramanujan, who independently obtained most of these classical results. Interpolation is very important also 
for the theory of zeta functions. Although for the solution of the Basel problem only one value at $s=2$ 
(Newton sum of squares) does matter, the interpolation to the real and the complex of the Newton sums, looking at the
exponent as parameter, is a fruitful \textit{transalgebraic} idea. Obviously the function $\zeta(2s)$ is the interpolation 
of the Newton sums of the roots of the $\frac{\sin(\pi \sqrt{z})}{\pi \sqrt{z}}$ transcendental function. 
Holomorphic functions of low order and with symmetries are uniquely determined by its values on arithmetic 
sequences of points as follows from 
Weierstrass theory of entire functions. Remarkably, under convexity assumptions, this also holds in the real domain 
(for example Bohr-Mollerup characterization of the $\Gamma$-function). It is precisely this point of view and 
the interpolation from the same points (real even integers) that leads to the $p$-adic version of the 
zeta function by Kubota and Leopold 
(see section \ref{sec:Kuwota-Leopold}).

\smallskip
 
There is another key point in Euler's argument. Multiplying 
$$
s(z)=\frac{\sin  \sqrt{z}}{ \sqrt{z}}=1-\frac{1}{3!} z+\frac{1}{5!} z^2 +\ldots
$$ 
by the exponential of an entire function changes the coefficients of the power series expansion at $0$, but 
not its divisor. We can keep these coefficients rational and the entire function of finite order by multiplication 
by the exponential of a polynomial with rational coefficients. Therefore we cannot expect that direct Newton relations 
hold for the roots an arbitrary entire function. But we know that this is true for entire functions of order $<1$ 
which can be 
treated as polynomials, as is the case for $s(z)$. Euler's transalgebraic point of view 
also ``proves" the result: $s(z)$ behaves as the minimal entire function for the transalgebraic
number $\pi^2$, similar to the minimal polynomial of an algebraic number in Algebraic Number Theory. Although this notion is 
not properly defined in a general context of entire function with rational coefficients, 
the idea is that it is the 
``simplest" (in an undetermined sense, but meaning at least of minimal order) entire function with rational coefficients 
having $\pi^2$ as root. Note here a major difference between even and odd powers of $\pi$, which is the origin 
of the radical computability difference of $\zeta (2n)$ and $\zeta (2n+1)$. This creative ambiguity between algebraic and transalgebraic permeates Euler's work as well as in some of his most brilliant successors 
like Cauchy and most notably Galois in his famous memoir on resolution of algebraic (transalgebraic?) equations \cite{Ga}. Only the transalgebaic aspects of Galois memoir can explain the failure of Poisson to understand some 
basic facts, and can explain some notations and claims, in particular of the lemma after the main theorem where 
Galois mentions the case of ``infinite number of roots''. 

Transalgebraicity is an extremely powerful ``philosophical" principle that some mathematicians of the XIXth century seem to be well aware of. 
In general terms we would say that 
{\it analytically unsound manipulations provide correct answers when they have an underlying
transalgebraic background}. This deserves to be vaguely stated as a general guiding principle:

\medskip

 \textbf{Transalgebraic Principle (TP).} {\it Divergent or analytically unsound algebraic manipulations yield correct results when they have a transalgebraic meaning.}

\medskip

We illustrate this principle with one of Euler's most famous example. 
The transalgebraic meaning of the exponential comes from Euler's classical 
formula: 

$$
e^z=\lim_{n\to +\infty} \left (1+\frac{z}{n} \right )^n \ .
$$
From Euler's viewpoint the exponential function is nothing else but a ``polynomial" with a single zero of 
infinite order at infinite. Thus we can symbolically write:

$$
e^z=\left (1+\frac{z}{\infty } \right )^{+\infty} \ ,
$$
noting that both infinities in the formula are distinct in nature: $\infty$ is a geometric point on the Riemann sphere (the zero) and $+\infty$ is the infinite number (the order).

Then we can recover the main property of the exponential by invoking to the Transalgebraic Principle and the following computation (using $\infty^{-2}<<\infty^{-1}$)

\begin{eqnarray*}
e^z .e^w &&=\left (1+\frac{z}{\infty } \right )^{+\infty} .\left (1+\frac{w}{\infty } \right )^{+\infty} \\
&&=\left (1+\frac{z+w}{\infty }+\frac{z.w}{\infty^2} \right )^{+\infty} \\
&&=\left (1+\frac{z+w}{\infty } \right )^{+\infty} \\
&&=e^{z+w} \ . \\
\end{eqnarray*}

\bigskip 
 
The next major appearance of the zeta function was in 1748, when Euler obtained his celebrated 
{\it product formula}, which is equivalent to  
the Fundamental Theorem of Arithmetic, as we see from the key manipulation,
\begin{eqnarray*}
\zeta(s)&&=\sum_{n=1}^\infty n^{-s}=\sum_{p_1,p_2,\ldots , p_r} \left (p_1\ldots p_r \right )^{-s}\\
&&=\prod_p \left (1+p^{-s}+p^{-2s}+\cdots \right )\\
&&=\prod_p \left (1-p^{-s}\right )^{-1}, 
\end{eqnarray*} 
where $p=2, 3, 5, 7, 11,\ldots$ runs over all prime numbers. 
The product formula provides the link between the zeta function and 
prime numbers. Or, if we prefer, once we consider the extension of Riemann zeta function to the complex plane, 
first proposed by Riemann, between prime numbers and complex analysis. The meromorphic extension of the zeta function 
to the complex plane, first proved by Riemann, is one of the main contributions to the theory and the reason why  
the name of Riemann is attached to the zeta function. Observe that 
complexification comes naturally once we realize that each local factor
\begin{equation*}
\zeta_p(s)=(1-p^{-s})^{-1}
\end{equation*}
has the \textbf{DL} property, since its divisor consists only of poles located at the purely imaginary axes 
$i \omega_p^{-1} \ZZ$, where 
$$
\omega_p=\frac{\log p}{2\pi} \ ,
$$
is the fundamental frequency associated to the prime $p$ (in other contexts $\omega_p^{-1}$ may be preferred as ''prime frequency").
Up to a proper complex linear rescaling (depending on $p$) and an exponential factor (which does not change the divisor), 
the local factor $\zeta_p^{-1}$ is nothing but the sine function studied before:
\begin{equation*}
\zeta_p(s)^{-1}=1-p^{-s}=2i p^{s/2} \sin \left(\pi \omega_p (-is)\right ) \ .
\end{equation*}
Observe also that from the very definition of prime number it follows that the frequencies $(\omega_p)_p$ 
are $\QQ$-independent. 

Already at this early stage, since the \textbf{DL} property is invariant by finite products, 
it seems natural to conjecture that the zeta function does also satisfy the \textbf{DL} property. 
We do not mean that the \textbf{DL} property should be invariant by infinite products in general. 
According to the Transalgebraic Principle \textbf{TP}, it is natural to conjecture this 
for \textit{Euler products} of arithmetical nature. There is little support, at this stage, 
for this conjecture. Later on more evidence will appear. Due to the complete lack of 
arguments in the literature in support for the Riemann Hypothesis, this weaker and related 
conjecture is worth mentioning for its appealing nature, and may well be one of the original 
motivations of B. Riemann for conjecturing the Riemann Hypothesis. We can formulate it in the following imprecise form:

\medskip

\textbf{Conjecture.}
\textit{The \textbf{DL} property is invariant by arithmetical Euler products.}

\medskip

This conjecture fits well from what we actually know of the scope of the Riemann Hypothesis: We only expect the 
Riemann hypothesis to hold when an Eulerian factorization exists. The reason is the Davenport-Heilbronn zeta function 
for which the Riemann Hypothesis fails, although it shares many properties with Riemann zeta function (see \cite{Iv2} section 3).
 
From the product formula, Euler derived a new analytic proof and strengthen Euclid's result on the infinity of prime 
numbers. We only need to work with the real zeta function. 
For real $s>1$ the product formula shows that $\zeta (s)$ remains positive, thus we can take its logarithm and compute by 
absolute convergence

\begin{equation*}
\log \zeta(s)=\sum_p -\log(1-p^{-s})=\sum_p \sum_{n=1}^\infty \frac{p^{-ns}}{n}=
\left (\sum_p p^{-s} \right )+\sum_{n=2}^\infty\frac{1}{n}\left (\sum_p p^{-ns} \right ) 
\end{equation*} 

Let $\cP (s)=\sum_p p^{-s}$ be the first series in the previous formula and $R(s)$ the remainder terms. 
Since $R(s)$ converges for $s>\frac{1}{2}$, and is uniformly bounded in a neighborhood of $1$, the divergence 
of the harmonic series $\zeta(1)=\lim_{s\rightarrow 1+}\zeta(s)=\sum_{n=1}^{+\infty}\frac{1}{n}$ is equivalent to that of 
$$
\cP(1)=\sum_p \frac{1}{p}=\infty \ .
$$ 
Thus, there are infinitely many prime numbers with a certain density, more precisely, enough to make divergent 
the sum of their inverses, as Euler writes:
\begin{equation*}
\frac{1}{1}+\frac{1}{2}+\frac{1}{3}+\frac{1}{5}+\frac{1}{7}+\frac{1}{11}+\ldots =+\infty \ .
\end{equation*}
Note that Euler includes $1$ as prime number.

\smallskip

\section{Riemann's article.}

As we have observed, complexification not only arises naturally from Euler work, but it is also necessary to justify his 
results. Euler's successful complexification of the exponential function proves that he was well aware of this fact. But the 
systematic study and finer properties of the zeta function in the complex domain appears first in Riemann's memoir on the 
distribution of prime numbers \cite{Ri}. Before Riemann, P.L. Tchebycheff \cite{Tc} used the real zeta function to prove rigorous bounds on 
the {\it Law of prime numbers}, conjectured by Legendre (and refined by Gauss): When $x\to +\infty$
\begin{equation*}
\pi(x)=\# \{p \ \ {\hbox{\rm prime }}; p\leq x \} \approx \frac{x}{\log x} \ .
\end{equation*}
Riemann fully exploits the meromorphic character of the extension of the zeta function to the complex plane in order to give an explicit formula for $\pi(x)$. We review in this 
section his foundational article \cite{Ri} that contain ideas and results that go well beyond the law of prime numbers.  

\subsection{Meromorphic extension.}

Riemann starts by proving the meromorphic extension of the $\zeta$-function from an integral formula.
From the definition of Euler $\Gamma$-function (we refer to the appendix for basic facts), for real $s>1$ we have
\begin{equation*}
\Gamma(s)\ \frac{1}{n^s}=\int_{0}^{+\infty}e^{-nx} x^s \ \frac{dx}{x} \ .
\end{equation*} 
From this Riemann derives
\begin{equation*}
\Gamma(s)\zeta(s)=\int_{0}^{+\infty}
\left (\sum_{n=1}^{+\infty}e^{-nx}\right ) \ x^s\ \frac{dx}{x}=\int_{0}^{+\infty}\frac{e^{-x}}
{1-e^{-x}} \ x^s \ \frac{dx}{x}=I(s) \ ,
\end{equation*} 
where the exchange of the integral and the series is justified by the local integrability of the 
function $x\mapsto \frac{x^{s-1}}{e^x-1}$ at zero. Therefore
\begin{equation*}
\zeta(s)=\frac{1}{\Gamma(s)}\int_{0}^{+\infty}\frac{1}{e^x-1}\ x^s \ \frac{dx}{x}=\frac{1}{\Gamma(s)} \ I(s) \ .
\end{equation*} 

\smallskip

In order to prove the meromorphic extension of $\zeta$ to the whole complex plane, Riemann observes that we can 
transform the above integral $I(s)$ into a Hankel type contour integral for which the integration makes sense for all
$s\in \CC$. For $z \in \CC-\RR_{+}$, and $s\in\CC$, we consider 
the branch $(-z)^s=e^{s\log(-z)}$ where the usual principal branch of $\log$ in $\CC-\RR_-$ is taken.

\begin{figure}[h] 
\centering
\resizebox{9cm}{!}{\includegraphics{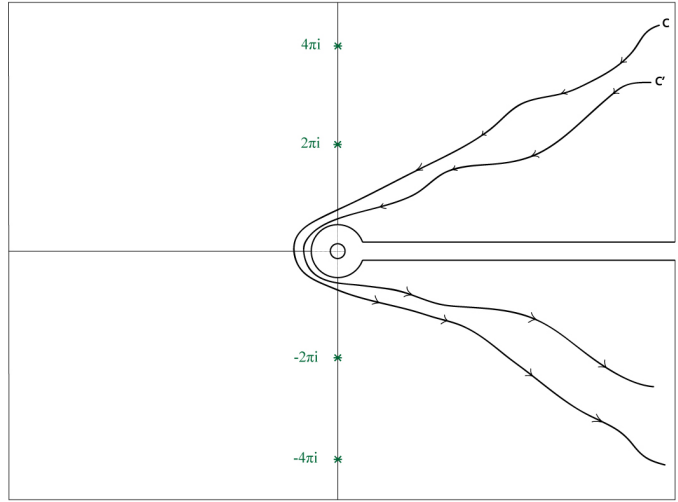}}    
\caption{Integration paths.}
\end{figure}

Let $C$ be a positively oriented curve as in Figure 1 uniformly 
away from the singularities $2\pi i \ZZ$ of the integrand, for example $d(C, 2\pi i \ZZ)>1$, with $\Re z \to +\infty$ 
at the ends of $C$, separating $0$ from $ 2\pi i \ZZ^*$, and with $0$ winding number for any singularity $2\pi i k$, 
with $k\not= 0$. Then the integral
\begin{equation*}
I_0(s)=\int_C \frac{(-z)^s}{e^z-1}\ \frac{dz}{z}
\end{equation*} 
defines an entire function which, by Cauchy formula, does not depend on the choice of the path $C$ in the 
prescribed class: Given another such path $C'$, we check that both integrals 
coincide by closing up both paths, in increasingly larger parts of $C$ and $C'$.

Now, for 
real $s>1$ and $\epsilon >0$ we can also deform the contour to obtain the expression
\begin{equation*}
I_0(s)=-\int_{\RR_{+}+i\epsilon} \frac{(-z)^s}{e^z-1}\ \frac{dz}{z}+\int_{|z|=\epsilon} \frac{(-z)^s}{e^z-1}\ \frac{dz}{z}+\int_{\RR_{+}-i\epsilon} \frac{(-z)^s}{e^z-1}\ \frac{dz}{z},
\end{equation*} 
which in the limit $\epsilon \rightarrow 0$ yields 
\begin{equation*}
I_0(s)=(e^{i\pi s}-e^{-i\pi s})I(s)=2i\sin(\pi s)I(s) \ ,
\end{equation*} 
where we have used that the middle integral tends to zero when $\epsilon\to 0$, because $s>1$ and for $z\to 0$
\begin{equation*}
\frac{(-z)^s}{e^z-1}=\cO (|z|^{s -1}) \ .
\end{equation*}
Therefore using the complement formula for the $\Gamma$-function,
\begin{equation*}
\zeta(s)=\frac{1}{2i\sin(\pi s)\Gamma(s)}\int_C \frac{(-z)^s}{e^z-1}\frac{dz}{z}=
\frac{\Gamma(1-s)}{2\pi i}\int_C \frac{(-z)^s}{e^z-1}\frac{dz}{z} \ .
\end{equation*} 
This identity has been established for real values of $s>1$, but remains valid by analytic continuation for all 
complex values $s\in \CC$. It results that $\zeta$ has a meromorphic extension to all of $\CC$ and its possible poles 
are located at $s=1,2,3,\ldots$. But from the definition for real values we know that $\zeta$ is positive and 
real analytic for real $s>1$, and becomes infinite at $s=1$, thus only $s=1$ is a pole. 
We can also see directly from the formula that  
$s=2,3,\ldots$ are not poles since the above integral is zero when $s\geq 2$ an integer 
(Hint: Observe that $(-z)^s$ has no monodromy 
around $0$ when $s$ is an integer, then use Cauchy formula). With the same argument, we get that $s=1$ is a pole 
with residue $\mathrm{Res}_{s=1}\zeta(s)=1$, since 
the residue of $\Gamma(1-s)$ at $s=1$ is $-1$ and by Cauchy formula
\begin{equation*}
I_0(1)=\int_C \frac{-z}{e^z-1}\frac{dz}{z}=-2\pi i  .
\end{equation*} 
We conclude that $\zeta(s)$ is a meromorphic function with a simple pole at $s=1$ with residue $1$.

It also follows from the integral formula that Riemann $\zeta$ function is a meromorphic function of order $1$, since 
the $\Gamma$-function is so and the integrand has the proper growth at infinite.

\subsection{Value at negative integers.}

To calculate the value of $\zeta$ at negative integers we introduce the Bernoulli numbers as the coefficients $B_n$ of the power series expansion
\begin{equation*}
\frac{z}{e^z-1}=\sum_{n=0}^{+\infty}\frac{B_n}{n!}z^n \ .
\end{equation*} 
Then we have $B_0=1$, $B_1=1$ and $B_{2n+1}=0$ for $n\geq 1$. The first values of Bernoulli numbers are
\begin{equation*}
B_2=\frac{1}{6}, \ \ B_4=-\frac{1}{30}, \ \ B_6=\frac{1}{42}, \ \ B_8=-\frac{1}{30}, \ \ B_{10}=\frac{5}{66}, \ \ B_{12}=-\frac{691}{2730} \ .
\end{equation*} 
As is often remarked, the appearance of the \textit{sexy} prime number 691 reveals the hidden presence of Bernoulli numbers. Making use of the above development, we can write
\begin{align*}
\zeta(1-n)&=\frac{\Gamma(n)}{2\pi i}\int_C \left( \sum_{k=0}^{+\infty}\frac{B_k}{k!}z^k \right) (-z)^{1-n} \ \frac{dz}{z^2} \\
&=(-1)^{1-n}\frac{(n-1)!}{2\pi i}\sum_{k=0}^{+\infty}\frac{B_k}{k!} \int_C z^{k-n-1}dz \\
&=(-1)^{1-n}\frac{B_n}{n} \ .
\end{align*} 
Thus $\zeta(1-n)=0$ for odd $n\geq 3$, and for $n$ even we have 
\begin{equation*}
\zeta(1-n)=-\frac{B_n}{n}\in \QQ \ .
\end{equation*}
In particular
\begin{equation*}
\zeta(0)=-\frac{1}{2} \ , \ \ \zeta(-1)=-\frac{1}{12} \ , \ \ \zeta(-3)=\frac{1}{120} \ , \ \ \zeta(-5)=-\frac{1}{252} \ , \ \ \zeta(-7)=\frac{1}{240} \ ,
\end{equation*} 
and $\zeta(-2k)=0$ for $k\geq 1$. These are referred in the literature as the \textit{trivial zeros} of $\zeta$. Note that these trivial zeros all lie on the real line and are compatible with the possible \textbf{DL} property for the $\zeta$ function.

\subsection{First proof of the functional equation.}

After establishing the meromorphic extension of the zeta function, Riemann proves a functional equation for $\zeta$ 
that is already present in the work of Euler, at least for real integer values. He gives two proofs. The first one that we 
present here is a direct proof from the previous integral formula. We assume first that $\Re s<0$.

\begin{figure}[h] 
\centering
\resizebox{8cm}{!}{\includegraphics{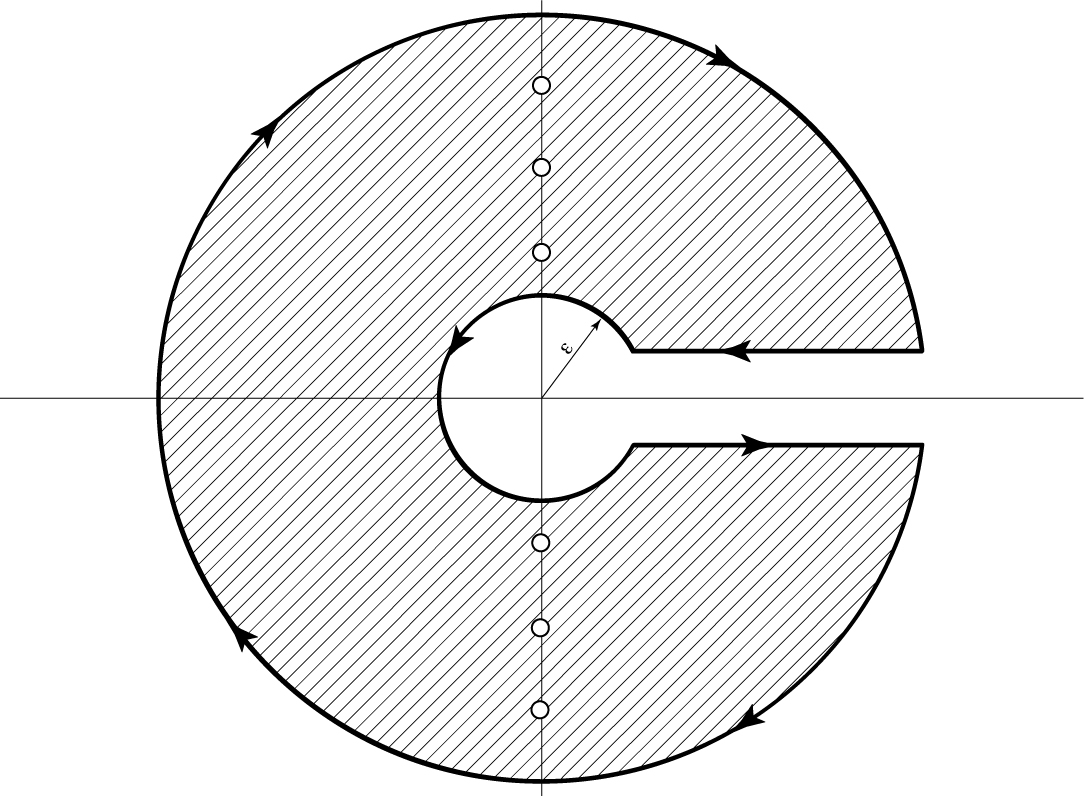}}    
\caption{Contour integral.}
\end{figure}

Applying the Residue 
Theorem to the contour integral in Figure 2, taking into account that the contribution of the outer boundary 
(second integral 
in the formula below) is zero in the limit (because $\Re s<0$), the residue theorem gives ($D_R$ is the outer 
circular part of the path and $C_R$ the rest):
\begin{align*}
I_0(s)&=\lim_{R\rightarrow +\infty}\left( \int_{C_R}\frac{(-z)^s}{e^z-1}\frac{dz}{z}+\int_{D_R}
\frac{(-z)^s}{e^z-1}\frac{dz}{z}\right) =-(2\pi i)\sum_{n\in \ZZ^*}\mathrm{Res}_{z=2\pi i n}\frac{(-z)^{s}}{e^z-1}\ \frac{1}{z} \\
&=-(2\pi i)\sum_{n=1}^{+\infty}\left( \frac{(2\pi n)^s e^{-\frac{s\pi}{2}i}}{2\pi i n}+\frac{(2\pi i n)^s e^{\frac{s\pi}{2}i}}{-2\pi i n}\right) =(2\pi)^s\sum_{n=0}^{+\infty}(e^{i\frac{s\pi}{2}}-e^{-i\frac{s\pi}{2}})\ n^{s-1} \\
&=(2\pi)^s 2i\sin\left (\frac{s\pi}{2}\right ) \zeta(1-s) \ .
\end{align*} 
Then $\zeta(s)=\frac{\Gamma(1-s)}{2\pi i}I_0(s)$ satisfies the functional equation,
\begin{equation}\label{eq:functional_equation} 
\zeta(s)=2^s\pi^{s-1}\sin\left (\frac{\pi s}{2}\right )\Gamma(1-s)\zeta(1-s) \ .
\end{equation} 
As an application, we can compute the values of $\zeta$ at even positive integers. For $s=1-2n,$ the functional equation gives
\begin{equation*}
\zeta(2n)=(-1)^{n+1}\frac{2^{2n-1}B_{2n}}{(2n)!}\pi^{2n} \ ,
\end{equation*} 
which is a rational multiple of $\pi^{2n}$. 

Equation (\ref{eq:functional_equation}) can be written in a more symmetric form, by using the division and complement formula for the gamma function:
\begin{align*}
\Gamma(1&-s)=2^{-s}\pi^{-\frac{1}{2}}\Gamma\left( \frac{1}{2}-\frac{s}{2}\right) \Gamma\left( 1-\frac{s}{2}\right) \ , \\
&\sin\left( \pi\frac{s}{2}\right) =\frac{\pi}{\Gamma\left( \frac{s}{2}\right) \Gamma\left( 1-\frac{s}{2}\right) } \ .
\end{align*} 
Then:
\begin{equation*}
\pi^{-\frac{s}{2}}\Gamma\left( \frac{s}{2}\right) \zeta(s)=\pi^{-\frac{1-s}{2}}\Gamma\left( \frac{1-s}{2}\right) \zeta(1-s) \ .
\end{equation*} 
One can multiply by $s(s-1)$ in order to kill the pole at $s=1$, respecting the symmetry, then 
if one introduces the extended zeta function
\begin{equation*}
\xi(s)=\frac{1}{2}s(s-1)\pi^{-\frac{s}{2}}\Gamma(s/2)\zeta(s) \ ,
\end{equation*} 
(the factor $1/2$ is only there for historical reasons because Riemann in his memoir was using the 
function $\Pi(s)=\Gamma (s-1)$), the functional equation becomes
\begin{equation*}
\xi(s)=\xi(1-s).
\end{equation*} 
The function $\xi$ is an entire function of order $1$ whose zeros are the zeros of $\zeta$ except for the trivial ones. 
Since $\xi(s)$ is non-zero when $\Re s>1$, the function does not vanish for $\Re s<0,$ so all its zeros belong to 
the \textit{critical strip} $0\leq \Re s \leq 1$. The symmetry of the equation can be formulated as the real analyticity of
$\Xi(t)=\xi\left( \frac{1}{2}+it\right) $, i.e. $\Xi(t)=\overline{\Xi(\overline{t})}$. Therefore we have two order two
symmetries. This makes that each non-real zero $\rho$ of $\zeta$ is associated 
to other zeros at the locations $1-\rho$, $\bar \rho$,
and $1-\bar \rho$. 

\subsection{Second proof of the functional equation.}

In his memoir, Riemann gives a second proof of the functional equation that discloses an important link with modular 
forms. This proof plays an important role in the theory of modular $L$-functions (see section \ref{sec:modular} and enjoy the reading of ``le jardin des d\'elices modulaires'', fourth volume of \cite{Go}). 
The functional equation indeed follows from the functional equation for the classical theta function.

We consider the classical $\theta$-function 
\begin{equation*}
\theta(x)=\sum_{n\in \mathbb{Z}} e^{-n^2\pi x^2} \ ,
\end{equation*}
and the half $\theta$-function
\begin{equation*}
\psi(x)=\sum_{n=1}^{+\infty} e^{-n^2\pi x^2} \ ,
\end{equation*}
so $\theta(x)=2\psi(x)+1$.

Again starting from the integral formula for the $\Gamma$-function, we have

\begin{equation*}
\pi^{-s/2}\Gamma(s/2)\frac{1}{n^s}=\int_0^{+\infty} e^{-n^2\pi x^2}x^{\frac {s}{2}} \ \frac{dx}{x} \ ,
\end{equation*}
and adding up for $n=1,2,3,\ldots$, 
\begin{equation*}
\pi^{-s/2} \Gamma(s/2)\zeta(s)=\int_0^{+\infty}\psi(x)x^{\frac{s}{2}}\ \frac{dx}{x} \ . 
\end{equation*} 

Recall Poisson formula: For a function $f:\RR\to \RR$ in the rapidly decreasing Schwartz class, we have
\begin{equation*}
\sum_{n\in \mathbb{Z}}f(n)=\sum_{n\in \mathbb{Z}} \hat{f}(n) \ .
\end{equation*} 
When we apply this to the Gauss function we obtain the classical functional equation 
$\theta(x)=x^{-\frac{1}{2}}\theta(1/x)$ or
\begin{equation*}
2\psi(x)+1=x^{-\frac{1}{2}}\left (2\psi(1/x)+1 \right )
\end{equation*} 
Thus splitting the integral in two parts, integrating over $[0,1]$ and $[1,+\infty [$, and using the functional equation 
$\psi(x) = x^{-1/2} \psi(1/x)+\frac{1}{2} x^{-1/2}-1/2$ we get 
\begin{equation*}
\pi^{-s/2} \Gamma(s/2)\zeta (s)=\int_1^{+\infty}\psi(x) x^{\frac{s}{2}}\ \frac{dx}{x} +\int_0^1 
\psi(1/x)x^{\frac{s-1}{2}} \ \frac{dx}{x}+ \frac{1}{2}\int_0^1 (x^{\frac{s-1}{2}}-x^{\frac{s}{2}})\ \frac{dx}{x} \ ,
\end{equation*}
then the change of variables $x\mapsto 1/x$ in the second integral gives the symmetric form, invariant by $s\mapsto 1-s$,

\begin{equation*}
\pi^{-s/2} \Gamma(s/2)\zeta (s)=-\frac{1}{s(1-s)}+\int_1^{+\infty}\psi(x)(x^{\frac{s}{2}}+x^{\frac{1-s}{2}})\ \frac{dx}{x} \ .
\end{equation*}
We can also write
\begin{equation*}
\xi (s)=\frac{1}{2}-\frac{s(1-s)}{2}\int_1^{+\infty}\psi(x)(x^{\frac{s}{2}}+x^{\frac{1-s}{2}})\ \frac{dx}{x} \ .
\end{equation*}
Then Riemann manipulates further the integral expression and writes $s=\frac{1}{2}+it$ and 
\begin{equation*}
s(1-s)=t^2+1/4 \ ,
\end{equation*}
\begin{equation*}
\Xi(t)=\xi\left( \frac{1}{2}+it\right) =\frac{1}{2}-2(t^2+1/4)^{-1}\int_1^{+\infty} \psi(x) x^{-3/4} \cos \left (\frac{1}{2} t \log x\right ) \ dx \ ,
\end{equation*}
which gives by integration by parts (and using the functional equation for $\psi$),
\begin{equation}\label{eq:xi_integral}
\Xi(t)=4\int_1^{+\infty} \frac{d \left (x^{-3/2}\psi(x) \right )}{dx} x^{-1/4}  \cos \left (\frac{1}{4} t \log x\right ) \ dx \ .
\end{equation}

\subsection{The Riemann Hypothesis.}

At that point of his memoir, Riemann observes that formula (\ref{eq:xi_integral}) shows that 
the function $\Xi$ can be expanded in a rapidly convergent power series of $t^2$. 
Then he makes the observation that the only zeros are in the strip $-i/2\leq \Im t \leq i/2$. 
Moreover he estimates the number 
of zeros of $\Xi$ with real part in $[0,T]$ explaining how to obtain the result by a contour integration 
and giving the asymptotic (Riemann- Von Mangoldt)
\begin{equation*}
 N(T) \approx \frac{T}{2\pi} \log \frac{T}{2\pi }-\frac{T}{2\pi } \ .
\end{equation*}
Then he continues making a curious observation:

\medskip

\textit{One now finds indeed approximately this number of real roots within these limits, and it is very likely 
that all roots are real.} 

\textit{Certainly, one would wish to have a rigorous proof of this proposition; but I have meanwhile temporarily 
left this question after some quick futile attempts, as it appears unnecessary for the immediate purpose of my investigations.}

\bigskip

In the manuscript notes that he left one can find the numerical computation of some of the first zeros 
(recall that $\g$ is a zero of $\Xi$ if and only if $\rho=1/2+i\gamma$ is a zero of $\xi$),
\begin{eqnarray*}
\gamma_1&&=14.134725\ldots, \\ 
\gamma_2&&=21.022039\ldots, \\ 
\gamma_3&&=25.01085 \ldots, \\ 
\gamma_4&&=30.42487 \ldots, \\
&&\vdots \\
\end{eqnarray*}

We refer to non-real zeros of $\zeta$ as \textit{Riemann zeros}.
The first part on Riemann's quotation seems to not have been properly appreciated: 
He claims to know that there are 
an infinite number of zeros for $\zeta$ in the critical line, and that 
if $N_0(T)$ is the number of zeros on the critical line with imaginary part in $[0,T]$, then when $T\to +\infty$
$$
N_0(T) \approx \frac{T}{2\pi} \log \frac{T}{2\pi }-\frac{T}{2\pi } \ .
$$
Only in 1914 G. H. Hardy proved that $N_0(T)\to +\infty$
when $T\to +\infty$, i.e. that there are an infinite number of zeros on the critical line. Hardy actually got the estimate,
$C>0$ universal constant 
\begin{equation*}
 N_0(T) > C \ T \ ,
\end{equation*}
which was improved in 1942 by A. Selberg to
\begin{equation*}
 N_0(T) > C \ T \log T \ ,
\end{equation*}
and later by N. Levinson (1974) to get an explicit $1/3$ proportion of zeros (and J.B. Conrey (1989) got $1/5$), but so 
far Riemann's claim that $N_0(T)\approx N(T)$ remains an open question. 

The second part in Riemann quotation attracted worldwide attention. He indicates that it is very likely that 
$N_0(T)=N(T)$ or in other terms, he is formulating what we now call the Riemann Hypothesis:

\bigskip

\textbf{Riemann Hypothesis}
\textit{Non-real zeros of the Riemann zeta function $\zeta$  are located in the critical line $\Re s=\frac{1}{2}$.}

\bigskip

Apart from his numerical study, we may ask what might be Riemann's intuition behind such claim. 
Obviously we can only 
speculate about this. As we see next, there seems to be some transalgebraic philosophical support, which may well
be the motivation behind his intuition.

We have seen that the functional equation forces the introduction of a natural complement 
$\Gamma$-factor with zeros in the real axes that satisfies the \textbf{DL} property (as well as the 
factor $s(1-s)$ which removes the pole at $s=1$, keeping the symmetry). Thus all the natural algebra 
related to Riemann zeta function involves functions in the \textbf{DL} class. This indicates that the whole theory
is natural in this class of functions.

The zeta function satisfies the \textbf{DL} property if its non-real zeros lie on vertical lines. 
These vertical lines should be symmetric with respect to $\Re s=1/2$. The strongest form of this property 
is when we have a minimal location of the non-real zeros in a single vertical line, that must necessarily be 
the critical line by symmetry 
of the functional equation. This is exactly what the  Riemann Hypothesis postulates.
Moreover, the transalgebraic principle \textbf{TP} coupled with the evidence that all Euler 
factors are in the \textbf{DL} class, gives support to the conjecture.

\medskip

Apparently the conjecture is called ``hypothesis" because Riemann's computations show that the location of the zeros give 
further estimates on the error term of the formula for the Density of Prime Numbers. Each non-trivial zero gives an oscillating term whose 
magnitude depends on the location of each zero, and more precisely on its real part. Under the Riemann Hypothesis one can show that for any $\epsilon >0$,
\begin{equation*}
 \pi(x) =\Li(x)+{\cO} (x^{\frac{1}{2}+\epsilon}) \ ,
\end{equation*}
where $\Li$ denotes the logarithmic integral
\begin{equation*}
 \Li(x)=\int_0^x \frac{dt}{\log t} =\lim_{\epsilon \to 0} \left ( \int_0^{1-\epsilon} \frac{dt}{\log t} + \int_{1+\epsilon}^x \frac{dt}{\log t} \right )\ ,
\end{equation*}
and we have the asymptotic development, when $x\to +\infty$,
\begin{equation*}
 \Li(x)\approx  \frac{x}{\log x} \sum_{n=0}^{+\infty} \frac{n!}{(\log x)^n}\approx \frac{x}{\log x} + \frac{x}{(\log x)^2}+\ldots 
\end{equation*}

Conversely, this estimate also implies the Riemann Hypothesis. Bounds on the error term and zero free regions 
in the critical strip are intimately related. 

The relation between Riemann zeta function and prime numbers is given by the Euler product, that unfortunately only 
holds in the half plane $\Re s >1$. A straightforward formula linking Riemann zeta function and the prime counting
function $\pi(x)$ is (see \cite{Ti} p.2), for $\Re s >1$,
\begin{equation*}
 \log \zeta(s)=s\int_2^{+\infty} \frac{\pi(x)}{x (x^s-1)} \ dx \ .
\end{equation*}

\medskip

The connection with the zeros can be readily seen by computing the logarithmic derivative, for real $s >1$

\begin{equation*}
 \log \zeta (s)=-\sum_p \log(1-p^{-s})=\sum_{n=1}^{+\infty} \sum_p \frac{1}{n} p^{-ns} \ ,
\end{equation*}
thus
\begin{equation*}
\frac{\zeta'(s)}{\zeta(s)}=-\sum_{p, n\geq 1} \log p \ p^{-ns} \ .
\end{equation*}
On the other hand

\begin{equation*}
\zeta(s)=\frac{\pi^{s/2}}{\Gamma(s/2)} . \frac{2}{s(s-1)} \ \xi(s) =\frac{\pi^{s/2}}{\Gamma(s/2)} . \frac{2}{s(s-1)} \ \xi(0) \ {\prod_\rho}^* \left( 1-\frac{s}{\rho}\right ) \ ,
\end{equation*}
where the product extends over the non-trivial zeros pairing symmetrical $\rho$ and $1-\rho$ to ensure convergence. 
The product expansion involving the zeros is only properly justified by Weierstrass and Hadamard work on entire functions,
but since Euler this type of expansions for entire functions of order $1$ is familiar to analysts. Disagreeing with other
authors, we have little doubt that Riemann was able to justify this point if needed be. 
In fact he mentions the asymptotic
density of zeros in order to justify the expansion and this is a key point. Computing 
the logarithmic derivative:
\begin{equation*}
\frac{\zeta'(s)}{\zeta(s)}=\frac{1}{2}\log \pi-\frac{1}{2}\frac{\Gamma'(s/2)}{\Gamma(s/2)}-\frac{1}{s}-
\frac{1}{s-1}+{\sum_\rho}^* \frac{1}{s-\rho} \ ,
\end{equation*}
and finally we get a relation that is the starting point for many results on the distribution of prime numbers:
\begin{align}\label{eq:prime-zero formula 1}
\sum_{n\geq 1} \Lambda(n) n^{-s}&= \sum_{p, n\geq 1} \log p \ p^{-ns} \notag \\
&=-\frac{1}{2}\log \pi+\frac{1}{2}\frac{\Gamma'(s/2)}{\Gamma(s/2)}+\frac{1}{s}+\frac{1}{s-1}-{\sum_\rho}^* \frac{1}{s-\rho} \ ,
\end{align}
where $\Lambda$ is the Von Mangoldt function defined by $\Lambda (n)=\log p$ if $n=p^m$ or $\Lambda (n)=0$ otherwise. 

\bigskip

Another indication that can give some support to the Riemann Hypothesis, and that was at the origin of 
famous unsuccessful attempts, is the formula
\begin{equation*}
 \frac{1}{\zeta(s)}=\prod_p \left ( 1-p^{-s} \right )=\sum_{n=1}^{+\infty } \frac{\mu (n)}{n^s} \ ,
\end{equation*}
where $\mu$ is M\"obius function: As follows from the formula, $\mu(n)=0$ if some square divides $n$, and $\mu (p_1 p_2\ldots p_l)=(-1)^l$,
for $p_1, p_2, \ldots  , p_l$ distinct primes.

Obviously the vanishing of $\mu$ at multiples of squares and its changing sign gives a better chance to have the series convergent for $\Re s >1/2$, something that obviously would imply the Riemann Hypothesis. This would equivalent to have the estimate for the Mertens function $M(x)$, for every $\epsilon >0$,
\begin{equation*}
 M(x)=\sum_{1\leq n\leq x} \mu(n) = \cO (x^{1/2+\epsilon}) \ .
\end{equation*}

Unfortunately, nobody has even been able to prove convergence for $\Re s> \sigma$ with $\sigma <1$.

\subsection{The Law of Prime Numbers.} The goal of Riemann's article is to present a formula for the law of prime numbers.
Starting from the formula, for real $s>1$, 
\begin{equation*}
\log \zeta (s)=-\sum_p \log(1-p^{-s})=\sum_{p, n\geq 1} \frac{1}{n} p^{-ns}  \ ,
\end{equation*} 
he gives an integral form to the sum on the right observing that
\begin{equation*}
p^{-ns}=s\int_{p^n}^{+\infty} \frac{dx}{x^{s+1}}=s\int_0^{+\infty} \mathbf{1}_{[p^n,+\infty)} \ x^{-s} \ \frac{dx}{x} \ ,
\end{equation*}
thus
\begin{equation}\label{eq:integral formula}
\frac{\log \zeta (s)}{s}=\int_1^{+\infty}\Pi(x)\ x^{-s}\ \frac{dx}{x} \ , 
\end{equation}
where $\Pi$ is the step function
\begin{equation*}
\Pi =\sum_{\genfrac{}{}{0pt}{}{p\geq 2}{ n\geq 1}} \frac{1}{n} \ \mathbf{1}_{[p^n,+\infty[} \ . 
\end{equation*}
Equation (\ref{eq:integral formula}) is the exact formula that Riemann writes, and, as noted in \cite{Fr}, this seems to indicate that Riemann, as Euler, considered the number $1$ as a prime number (but the sums on primes do start at $p=2$).
The jumps of $\Pi$ are positive, located at the powers $p^n$ of primes and have magnitude $1/n$, thus we can also write
\begin{equation*}
\Pi(x)=\pi(x)+\frac{1}{2}\pi(x^{\frac{1}{2}})+\frac{1}{3}\pi(x^{\frac{1}{3}})+\cdots =\sum_{n=1}^{+\infty }  \frac{1}{n} \pi (x^{\frac{1}{n}}) \ ,
\end{equation*} 
where as usual $\pi(x)$ is the number of primes less than or equal to $x$. Note that $\pi(x)$ can be recovered from $\Pi(x)$ 
by M\"obius inversion formula:
\begin{equation*}
\pi(x)=\sum_{n=1}^{+\infty }  \frac{\mu(n)}{n} \Pi (x^{\frac{1}{n}})=\Pi(x)-\frac{1}{2}\Pi(x^{\frac{1}{2}})-\frac{1}{3}\Pi(x^{\frac{1}{3}})+\frac{1}{6}\Pi(x^{\frac{1}{6}})+\cdots  \ .
\end{equation*}

Now the integral in \ref{eq:integral formula} is a Mellin transform, or a Fourier-Laplace transform in the $\log x$ variable
\begin{equation}
\frac{\log \zeta (s)}{s}=\int_1^{+\infty}\Pi(x)\ e^{-s \log x}\ d (\log x) \ , 
\end{equation}
and Riemann appeals to the Fourier inversion formula to get, previous proper redefinition of $\Pi$ at the jumps (as the average of left and right limit as usual in Fourier analysis), $\Pi^*(x)=(\Pi(x+0)+\Pi(x-0))/2$,
\begin{equation}\label{eq:inver}
\Pi^*(x)=\frac{1}{2\pi i}\int_{a-i\infty}^{a+i\infty } \frac{\log \zeta (s)}{s}x^s \ ds \ ,
\end{equation} 
where the integration is taken on a vertical line of real part $a>1$. 

Next, taking the logarithm of the identity
\begin{equation*}
\zeta(s)=\frac{\pi^{s/2}}{\Gamma(s/2)} . \frac{2}{s(s-1)} \ \xi(s) =\frac{\pi^{s/2}}{\Gamma(s/2)} . \frac{2}{s(s-1)} \ \xi(0) \ {\prod_\rho}^* \left( 1-\frac{s}{\rho}\right ) \ ,
\end{equation*}
we get
\begin{equation*}
\log \zeta(s)=\frac{s}{2}\log\pi -\log(\Gamma(s/2))-\log (s(s-1)/2)+\log \xi(0)+{\sum_\rho}^* \log \left (1-\frac{s}{\rho} \right ) \ .
\end{equation*}
We would like to plug this expression into the inversion formula. Unfortunately there are convergence problems. Riemann uses the classical trick in Fourier analysis of making a preliminary integration by parts in order to derive from (\ref{eq:inver}) a convenient formula where the previous expression can be plug in, and the computation carried out integrating term by term (but always pairing associated non-trivial zeros),
\begin{equation}
\Pi^*(x)=-\frac{1}{2\pi \log x} \int_{a-i\infty}^{a+i\infty } \frac{d}{ds}\left (\frac{\log \zeta (s)}{s}\right ) 
\ x^s \ ds \ .
\end{equation} 
Once we replace $\log \zeta(s)$ and we compute several integrals, we get (see \cite{Fr} \cite{Ed} for details) the main result in Riemann memoir:
\begin{equation*}
\Pi^*(x)=\Li(x)-\sum_{\Im \rho >0} \left ( \Li (x^\rho)+ \Li (x^{1-\rho}) \right ) +\int_x^{+\infty}\frac{dt}{t(t^2-1)\log t}-\log(2) \ ,
\end{equation*} 
thus if the oscillating part coming from the non-trivial zeros can be neglected, the logarithmic integral 
\begin{equation*}
\Li(x)=\int_0^x \frac{dx}{\log x} 
\end{equation*} 
gives the principal part. Note that the integral on the right side converges to $0$ when $x\to +\infty$ and thus is irrelevant for the asymptotic. Note also that the series cannot be uniformly convergent since the sum is a discontinuous function. 
Finally, by M\"obius inversion we get, denoting $\pi^*(x)=(\pi(x+0)+\pi(x-0))/2$,
\begin{equation*}\label{eq:asymptotic}
\pi^*(x)=\sum_{n=1}^{+\infty } \frac{\mu(n)}{n}\ \Li(x^{\frac{1}{n}}) +R(x) \ ,
\end{equation*} 
where the remainder $R(x)$ is given, up to bounded terms, by the oscillating part coming from the non-trivial zeros
\begin{equation*}
R(x)=\sum_{n=1}^{+\infty } \sum_\rho \frac{\mu(n)}{n} \ \Li(x^{\frac{\rho}{n}}) +\cO (1) \ ,
\end{equation*} 
The first term in (\ref{eq:asymptotic}) suggest the asymptotic when $x\to +\infty$,
\begin{equation*}
 \pi (x) \approx \Li (x) \approx \frac{x}{\log x} \ ,
\end{equation*}
but Riemann didn't succeed in carrying out a control of the remainder, and the result was only proved rigorously independently by J. Hadamard and C.J. de la Vall\'e Poussin in 1896 (see \cite{Ha} and \cite{VP1}, \cite{VP2}), thus settling the long standing conjecture of Legendre and Gauss. The methods of Hadamard and de la Vall\'ee Poussin give the estimate,

\begin{equation*}
 \pi(x)=\Li(x) +\cO \left (x e^{-c \log x} \right ) \ , 
\end{equation*}
where $c>0$ is a constant. The Riemann Hypothesis would be equivalent to the much stronger estimate
\begin{equation*}
 \pi(x)=\Li(x) +\cO (\sqrt{x} \log x ) \ . 
\end{equation*}

Only more than half a century later, in 1949, P. Erd\"os and A. Selberg succeeded in giving an proof that didn't 
require the use of complex analysis. But the best bounds on the error do require the use of complex variable methods.

\medskip

Nevertheless, a more precise asymptotic development at infinite is suggested by Riemann's results. It is interesting to note that the main purpose of the first letter that S. Ramanujan wrote to G.H. Hardy was to communicate his independent discovery of the asymptotic (see \cite{Ra})
\begin{equation*}
\pi(x)\approx \sum_{n=1}^{+\infty } \frac{\mu(n)}{n}\ \Li(x^{\frac{1}{n}}) \ ,
\end{equation*} 
that Ramanujan claimed to be accurate up to the nearest integer, something that numerically is false. His exact words are the following (\cite{Har} p.22):

\medskip

\textit{ I have found a function which exactly represents the number of prime numbers less than $x$ , ``exactly" in the sense that the difference between the function and the actual number of primes is 
\textbf{generally} $0$ or some small finite value even when $x$ becomes infinite...(Letter from Ramanujan to Hardy)}

\medskip

Hardy, much later, after Ramanujan passed away, make frivolous comments about Ramanujan claim, and took 
good care of explaining why he believed that Ramanujan didn't seem to know about complex zeros of 
Riemann $\zeta$-function (\cite{Har} chapter II). 
Considering Ramanujan success with similar claims, as for instance for the asymptotic of partition numbers (see \cite{Rad}),
one should probably not take the words of Ramanujan so lightly. It is reasonable to conjecture,
interpreting his ''generally'', that probably Ramanujan claim
holds for arbitrarily large values of $x\in \RR_+$ :

\begin{conjecture}
 Let $\cR\subset \RR_+$ be the Ramanujan set of values $x\in \RR_+$ for which 
 $\pi(x)$ is the nearest integer to 
 \begin{equation*}
  \sum_{n=1}^{+\infty} \frac{\mu(n)}{n}\ \Li(x^{\frac{1}{n}}) \ ,
 \end{equation*}
 then $\cR$ contains arbitrarily large numbers. 
\end{conjecture}

In general one may ask how large is the Ramanujan set $\cR$. Does it have full density? That is:
 \begin{equation*}
  \lim_{X\to +\infty}\frac{\left |\cR \cap [0,X]\right |}{X} = 1 \ .
 \end{equation*}

This does not seem to agree with the plot in \cite{Za}, but one can expect to have 
elements of the Ramanujan set in large gaps of primes (large according to how far away is the gap).


\section{On Riemann zeta-function after Riemann.}

A vast amount of work has been done since Riemann's memoir. It is out of the scope of these 
notes to give an exhaustive survey. 
The aim of this section is to describe some results that we consider the most
relevant ones that shed some light or might be important to the resolution of the Riemann Hypothesis.
In particular, we will not discuss all the numerous equivalent reformulations of the Riemann Hypothesis.
Some of these are truly amazing, but unfortunately none of them seem to provide any useful insight 
into the conjecture.

\subsection{Explicit formulas.}

Riemann's formula for $\pi(x)$ in terms of expression involving the non-trivial zeros of 
$\zeta$ is an example of explicit formula for the zeros. It is surprising that despite how 
little we know about the exact location of the zeros, many such formulas do exist.

Already Riemann stated the formula, for $x>1$ and $x\not= p^m$,
\begin{equation*}
 \sum_{n\leq x} \Lambda (n)=x-\sum_\rho \frac{x^\rho}{\rho} -\frac{\zeta'(0)}{\zeta(0)}-\frac{1}{2}\log(1-x^{-2}) \ ,
\end{equation*}
where the sum over the zeros is understood in the sense
\begin{equation*}
\sum_\rho \frac{x^\rho}{\rho} =\lim_{T\to +\infty} \sum_{|\Im \rho| \leq T} \frac{x^\rho}{\rho} \ .
\end{equation*}
This formula can be proved applying the integral operator 
\begin{equation*}
 f\mapsto \frac{1}{2\pi} \int_{a-i\infty}^{a+i\infty} f(s) \ x^s \ \frac{ds}{s} \ .
\end{equation*}
to formula (\ref{eq:prime-zero formula 1}). 

Another very interesting explicit formula was given by E. Landau (\cite{La2} 1911), that shows 
the very subtle location of the zeros and how they exhibit resonance phenomena at the prime frequencies:
When $x>1$ is a power of $p$, $x=p^m$, we have,
\begin{equation*}
 \sum_{0<\Im \rho <T} x^{\rho} = -\omega_p \ T + \cO (\log T) \ ,
\end{equation*}
but for $x\not=p^m$ for any prime, we have,
\begin{equation*}
 \sum_{0<\Im \rho <T} x^{\rho} = \cO (\log T) \ .
\end{equation*}
Notice that in the first case we have an equivalent which is $\cO(T)$ and 
the terms of the sum do not cancel each other at the same rate as in the second
case where the equivalent is $\cO(\log T)$. 
Related to this sum,  H. Cramer \cite{Cr} studied the analytic 
continuation of the now called Cramer function
\begin{equation*}
 V(t)=\sum_{\Im \rho >0} e^{\rho t} \ .
\end{equation*}
The series is convergent in the upper half plane and despite the singularity at $0$ and simple poles at 
$\ZZ^* \log p$, Cramer proved its meromorphic continuation. He determined the monodromy at $0$ and its poles and residues. A.P. Guinand \cite{Gu2} proved the Guinand's functional equation which is at the source of the meromorphic
extension across the real axes:
\begin{equation*}
 e^{-t/2} V(t)+e^{t/2} V(-t)=e^{t/2}+e^{-t/2}-e^{-t/2} \ \frac{1}{1-e^{-2t}} \ .
\end{equation*}

\bigskip

The general technique to obtain these explicit formulas and results is   
by computing a contour integral $\int_R \zeta'(s) /\zeta (s) \ f(s) \ ds$ 
over a rectangular contour $R$, symmetric with respect to the critical line,  
with horizontal and vertical sides away from the critical strip in two ways: Using the residue 
theorem (the integrand has simple poles at the non-trivial zeros), and by integrating over the vertical sides,
exploiting the symmetry given by the functional equation and the Euler product expansion on the 
right side (the integral over the horizontal sides is shown to be negligible or tending to $0$). 
This is also the technique used to count the zeros in the 
critical strip with bounded imaginary part. One can find many of these explicit formulas in the two volume 
treatise of Landau \cite{La}, or in Ingham's book \cite{In}.

\medskip

First A.P. Guinand \cite{Gu1} and then J. Delsarte \cite{De}, noticed the Fourier-Poisson duality that 
transpire from these formulas. For a modern discussion on this aspect and the relation to trace formulas, 
we refer to P. Cartier
and A. Voros enlightening article \cite{CV}. Similar to Poisson formula that can be formulated as a 
distributional result about the Fourier transform of a Dirac comb, one can give a general distributional 
formulation, and the application to various test functions yield the different classical 
\textit{explicit formulas} with the Riemann zeros.

\begin{theorem} \textbf{(Delsarte Explicit Formula)} \label{thm_delsartre} In the distribution sense, 
 \begin{equation*}
 2\pi \left ( \delta_{i/2} + \delta_{-i/2} + \sum_\gamma \delta_\gamma \right )=-\sum_{n=1}^{+\infty} \frac{\Lambda (n)}{\sqrt{n}} \ \left ( \hat \delta_{\frac{\log n}{2\pi}} + \hat \delta_{-\frac{\log n}{2\pi}} \right ) - D_\infty  \ ,
\end{equation*}
where $D_\infty$ is the absolutely continuous distribution
\begin{equation*}
 D_\infty (\varphi )=\int_\RR \varphi (t) \Psi(t) \ dt \ ,
\end{equation*}
where 
\begin{equation*}
 \Psi(t)= \log \pi -\frac{1}{2} \frac{\Gamma'(1/2+it)}{\Gamma(1/2+it)} -\frac{1}{2} \frac{\Gamma'(1/2-it)}{\Gamma(1/2-it)} \ .
\end{equation*}
\end{theorem}

In 1952 A. Weil gave a general adelic version of this explicit formula. The form of the distribution corresponding to 
the infinite place was simplified by K. Barner (see \cite{We} and \cite{Lan}). Weil observed that the Riemann 
Hypothesis is equivalent to the positivity of the right hand side of his equation for a suitable class of functions (see \cite{Bo}). 
This type of positivity criteria is common to other approaches to the Riemann Hypothesis, as for instance De Branges proposed approach \cite{DB}.

\textit{Explicit formulas} have been interpreted as \textit{trace formulas} and conversely. The 
\textit{trace formula} approach was inaugurated by Selberg's work \cite{Se} (1956) on the length spectrum
of constant negatively curved compact surfaces. It is usually presented as a relation between the Fourier transform of 
distribution associated to the length of primitive geodesics, and 
the eigenvalues of the laplacian. His motivation is by analogy with explicit formulas for the Riemann zeta function.
Cartier and Voros gave an interpretation as an explicit formula in 
\cite{CV}. All of this is very much related to the procedure of zeta-regularization that is 
cherished by physicists (\cite{El}, \cite{EORBZ}), number theorists (\cite{Den}), an first used by geometers (\cite{RS})
in order to compute the determinant of the laplacian in a surface with constant negative curvature.

Another geometric interpretation of explicit formulas as a \textit{trace formula} in a non-commutative space
was proposed by A. Connes \cite{Con}. J.-B. Bost and A. Connes proposed 
such a space in \cite{BC}. Connes reduced the Riemann Hypothesis to the proof of a trace formula. 
This is very much related to \textit{quantum chaos} (see section \ref{sec:chaos}).

Recently, Selberg's analogy of his trace formula with explicit formulas for the Riemann zeta funciton was 
elucidated by V. Mu\~noz and the author \cite{MPM}. For any meromorphic Dirichlet series there is a general 
distributional formula that relates the zeros and poles
to the frequencies of the Dirichlet series: The Poisson-Newton formula. This formula unifies in one formula Poisson formula 
in Fourier analysis and classical Newton relations between power sums of roots and elementary symmetric functions: It is 
fundamentally transalgebraic in its nature. This corroborates the accuracy of the transalgebraic point of view in the
Riemann zeta function theory.

\subsection{Poisson-Newton formula}\label{sec:poisson-newton}

We consider a Dirichlet series $f(s)=1+\sum_{n\geq 1} a_n e^{-\lambda_n s}$ absolutely convergent 
in a right half plane, with $a_n\in \CC$, $0<\lambda_1 <\lambda_2<\ldots$
where the set of frequencies $(\lambda_n)$ are finite or $\lambda_n \to +\infty$. We assume that $f$ has a meromorphic 
extension of finite order to the whole complex plane. The distribution
$$
W(f)=\sum_\rho n_\rho e^{\rho t}
$$
where the sum runs over the divisor of $f$, with $n_\rho \in \ZZ$ being the multiplicity of $\rho$,
is a well defined in $\RR_+^*$ and can be naturally extended to a distribution in $\RR$. We call it 
the \emph{Newton-Cramer distribution} of $f$.

We can define the coefficients $(b_{\bk})$ by 
\begin{equation} \label{eqn:bn}
-\log f(s)=-\log \left ( 1+ \sum_{n\geq 1} a_n \ e^{-\lambda_n s}\right )
=\sum_{\bk \in \Lambda} b_{\bk} \, e^{-\langle \boldsymbol{\lambda} , \bk \rangle s}
 \ ,
 \end{equation}
where $\Lambda=\{ \bk=(k_n)_{n\geq 1} \, | \, k_n \in \NN, ||\bk||=\sum | k_n |<\infty, ||\bk|| \geq 1\}$,
and
$\langle \boldsymbol{\lambda} , \bk \rangle = \lambda_1k_1+\ldots + \lambda_{l}k_{l}$, where
$k_n=0$ for $n>l$. This is well defined for $\Re s \to +\infty$.

\begin{theorem}[Poisson-Newton Formula]
We have
$$
W (f)= \sum_{k=0}^{g-1} c_k 
 \delta_0^{(k)} + \sum_{\bk \in \Lambda} \langle \lambda , \bk
 \rangle \, b_{\bk} \ \delta_{\langle \boldsymbol{\lambda} ,\bk\rangle } \, .
$$
\end{theorem}

The structure of the distribution at $0$ is related to the exponential factor in Hadamard-Weierstrass factorization of $f$.

The amazing property of this formula is its duality. Applying it to the simpler Dirichlet series $f(s)=1-e^{-\lambda s}$ it 
gives the classical distributional form of Poisson formula in Fourier analysis,

\begin{equation}\label{eqn:a1}
\sum_{n\in \ZZ} e^{i\frac{2\pi}{\lambda} n t}  =\lambda \sum_{k\in \ZZ}  \delta_{\lambda k} \, ,
\end{equation}

When we apply the Poisson-Newton formula to the finite Dirichlet series with one frequency
$$
f(s)=e^{-\lambda ns }P(e^{\lambda s}) =1+a_1 e^{-\lambda s} +\ldots + a_n e^{-\lambda n s} \ .
$$
where $P(z)=z^n+a_1 z^{n-1}+\ldots +a_n $ is a polynomial, we 
obtain all Newton relations between the roots and the coefficients of $P$.
This transalgebraic duality somewhat explains the intuitions of Delsartre, Selberg, Cartier, Voros, etc

When we apply Poisson-Newton formula to $f(s)=\zeta(s)$ we obtain Delsartre's Explicit Formula (Theorem \ref{thm_delsartre}). 
When we apply Poisson-Newton formula to Selberg zeta function we get Selberg trace formula (see section \ref{sec:dynamical}), 
which proves that the origin of both formulas is the same, confirming Selberg's classical analogy.

Despite that popular proofs of the Explicit Formula uses the functional equation of Riemann $\zeta$-function, it is 
an interesting  feature that Poisson-Newton formula holds independently of having a functional equation.

Also the Poisson-Newton formula is indeed more general than stated. We can apply it to meromorphic functions $f$ of 
finite order with divisor contained in a half plane that are not Dirichlet series. The general theorem is the following:

\begin{theorem}[General Poisson-Newton Formula]
We have
$$
W(f)=\sum_{j=0}^{g-1} c_j \delta_0^{(j)} + \cL^{-1} (f'/f) \ ,
$$
where $\cL^{-1}$ is the inverse Laplace transform 
\end{theorem}

We can apply this general formula to the $\Gamma$ function for example. The result gives the classical Gauss integral formula for 
the logarithmic derivative of the $\Gamma$ function:

$$
\frac{\Gamma'(s)}{\Gamma (s)}=\int_0^{+\infty} \left (\frac{e^{-t}}{t}-\frac{e^{-st}}{1-e^{-t}}
\right ) \, dt \ .
$$
So, amazingly, we can interpret Gauss integral formula as all Newton relations 
for the roots $e^{-n}$ for all exponents $t\in \RR_+^*$. So this is the form that Newton relations can have in the transalgebraic context.

The analysis of the structure at $0$ of the Newton-Cramer distribution gives also an infinite number of new 
McLaurin type formulas and is very much related to Ramanujan treatment of infinite series (Ramanujan's 
theory of the constant).

\subsection{Montgomery phenomenon.}\label{sec:Montgomery}

In 1972 H. Montgomery \cite{Mo} studied the distribution of differences of Riemann zeros, and proved, assuming the Riemann Hypothesis plus other assumptions, a weak version of the \textit{Pair Correlation Conjecture} for the zeros 
$\rho=1/2+i \gamma$, for $T\to +\infty$, 
\begin{align}\label{eq:GUE}
 &\frac{1}{N(T)} \# \{ (\gamma , \gamma' ) ; 0<\gamma, \gamma' <T \ ; \gamma\not=\gamma' ; \ \frac{2\pi }{\log T} \alpha < \gamma -\gamma' < \frac{2\pi }{\log T} \beta \} \notag 
 \\ &\to \int_\alpha^\beta \left (1-\left (\frac{\sin \pi x}{\pi x} \right )^2 \right ) \ dx \ .
\end{align}

One finds the claim in the literature is that this statistic behaviour 
gives support to the so called Hilbert-Polya approach to the Riemann Hypothesis: It is 
enough to show that the Riemann zeros form the spectrum of an Hermitian operator. For the Hilbert-Polya
attribution one may consult the correspondence of Odlyzko with Polya in \cite{Od2}. It seems that E. Landau 
asked Polya for a physical reason that would explain the Riemann Hypothesis and he came out with his proposal
of a physical system whose eigenvalues (energy levels) correspond to the Riemann zeros.
As the story tells, F. Dyson recognized the same distribution as for the pair correlation of 
eigenvalues for the Gaussian Unitary 
Ensemble (GUE) in Random Matrix Theory used by physicists, which is the set of 
Hermitian matrices endowed with a gaussian measure.

Extensive computations (O. Bohigas, M. Gianonni \cite{BG}, and A. Odlyzko \cite{Od} \cite{Od2}) have 
confirmed the statistics with high accuracy. In particular, A. Odlyzko has pushed the computation of 
Riemann zeros to various millions and imaginary parts of the order of $10^{20}$. 

In the personal opinion of the author, it is a bold statement to claim that the GUE statistics (\ref{eq:GUE}) 
infers anything about Hilbert-Polya proposal. GUE statistics may just mean that 
the same Fourier phenomena arises in Random Matrix 
Theory and other domains. The GUE distribution differs from the constant distribution by 
the square of a Fresnel or sine-cardinal (non-positive) distribution 
\begin{equation}\label{eq:Fresnel}
\int_\alpha^\beta \frac{\sin(\pi x)}{\pi x} \ dx \ , 
\end{equation}
which is just the Fourier transform of a box (or truncation) function. The self-convolution of a box function 
is a tent function and this is essentially what we can read in the GUE distribution.

In view of the existing literature, it seems to have been unnoticed that the properly 
normalized error term in (\ref{eq:GUE})  is asymptotic when $T\to +\infty$ to a superposition of 
Fresnel distributions \cite{PM0}, one positive with larger amplitude and the substraction of an infinite series of 
lower order ones.

\subsection{Quantum chaos.}\label{sec:chaos}

Polya's proposition to Landau to seek for a physical system having the Riemann zeros as energy levels, 
was pushed forward by M. Berry (\cite{Be}) and then by S. Keating (\cite{Ke}, \cite{BK}). They assume the 
Riemann Hypothesis, all the approach is perturbative, and carry out liberally the computations 
using the divergent Euler product in the critical strip. The purpose is to explore 
the semiclassical formalism, 
and how it fits with the classical formulas for the Riemann zeta-function. 

The approach consists in considering the imaginary part of the Riemann zeros as energy levels 
of a semi-classical system and 
exploit classical perturbative techniques for the quantification of hyperbolic dynamical systems 
(this explains the term ``quantum chaos''). Indeed
there is now a well developed theory of quantification of chaotic dynamical systems (see \cite{Guz}) which
was inspired from classical physical techniques (as WKB method). A main result is Gutzwiller Trace Formula that 
relates the energy levels of the quantification of the classical system with the classical periodic orbits. Only 
for specific dynamical systems with simple dynamical $\zeta$-functions (for example for hyperbolic systems, 
see section \ref{sec:dynamical}) a rigorous formula
can be established.
The starting point is the observation that
$N(E)$ can be seen as the integrated density of states and Riemann-Von Mangoldt formula gives the asymptotic
\begin{equation*}
 N(E)=\frac{E}{2\pi}\log \left (\frac{E}{2\pi}\right ) -\frac{E}{2\pi} +\frac{7}{8}+ \frac{1}{\pi} \Im \log \zeta(1/2+iE) +\cO(1) \ .
\end{equation*}
This expression is composed by the averaged part
\begin{equation*}
 < N(E)>=\frac{E}{2\pi}\log \left (\frac{E}{2\pi}\right ) -\frac{E}{2\pi} +\frac{7}{8} \ ,
\end{equation*}
and the oscillatory term
\begin{equation*}
 N_{osc}(E)=\frac{1}{\pi} \Im \log \zeta(1/2+iE) \ .
\end{equation*}
Plugging in brutally Euler product, replacing $\zeta(1/2+iE)$, we get
\begin{equation*}
 N_{osc}(E)=-\frac{1}{\pi} \sum_p \sum_{k\geq 1} \frac{1}{k} e^{-\frac{1}{2} k \log p} \sin(E k \log  p) \ .
\end{equation*}
The density of states distribution is obtained by differentiating 
\begin{equation*}
 d(E)=dN(E)/dE=\bar d(E) -2 \sum_p \sum_{k\geq 1} \frac{\log p}{2\pi}e^{-\frac{1}{2} k \log p} \cos(E k \log  p)
\end{equation*}
In the semiclassical approximation of a hamiltonian system one can also give an equivalent when 
$\hbar \to 0$ to the density of states and we get an averaged part and an oscillatory part that 
can be computed using the periodic 
orbits of the system and its Lyapunov exponents:
\begin{equation*}
 d(E)=\bar d(E)+\frac{2}{\hbar} \sum_p \sum_{k\geq 1} A_{p,k}(E) \cos\left (\frac{k}{\hbar} 
 S_p(E)-k\frac{1}{2} \alpha_p \right ) \ ,
\end{equation*}
where $k$ is number of times the primitive orbit is counted, 
the phase $\alpha_p$ is the Maslov index around the orbit (monodromy in the symplectic linear group of the linearization along the orbit), 
and $S_p(E)$ is the action around
the orbit 
\begin{equation*}
 S_p(E) =\oint {\bf p} . d {\bf q} \ .
\end{equation*}
The amplitudes are given by
\begin{equation*}
 A_{p,k} =\frac{T_p(E)}{2\pi} e^{-\frac{1}{2} \lambda_p(E) k T_p(E)} \ ,
\end{equation*}
where $T_p(E)$ is the period of the primitive orbit, and $\lambda_p$ is the sum of Lyapunov exponents (trace 
of the linear part of the Poincar\'e return map). 

\medskip

This allows to infer some information about the hypothetical Polya hamiltonian system. The periods of the 
primitive orbits would be the logarithms of primes $T_p=\log p$, $\alpha_p=0$ and $\lambda_p=1$.
Unfortunately, there are also notable differences, the most important one the divergent sign of the oscillatory
part. But these computations have proved useful and serve as a guide to infer results, 
as for example the prediction of $n$-point correlations between Riemann zeros.

\medskip

The amazing analogy means one of two things: Or there is really a Polya physical system behind Riemann zeta function, 
or both theories share a good deal of common formalism. Nobody has been able to have a good 
guess for such a Polya system for the $\zeta$-function, and if it exists it would be very fundamental. One would expect  
Polya systems in other cases where the Riemann Hypothesis is expected, as for example for Dirichlet $L$-functions (section \ref{sec:Dirichlet}). The Polya 
system for the $\zeta$-function would be the most basic. 
We believe that Poisson-Newton Formula (section \ref{sec:poisson-newton}) has shed some light on this question.
It indicates we have a more general theory that all this explicit or trace formulas at ounce, and that there is a general theory for some aspects above
dynamical or physical $\zeta$-functions and arithmetic $\zeta$-functions. So, for instance, Poisson-Newton formula gives Guntzwiller Trace Formula when we 
have a dynamical $\zeta$ function with good properties and 
the Explicit Formula for Riemann $\zeta$ function at the same time. In view of this and the pointed out mismatchs,  
we are inclined to favor the second possibility, and that a natural Polya system for the Riemann $\zeta$-function may not exist.

\subsection{Statistics on Riemann zeros.}

We describe in this section some numerical results of transalgebraic origin from the author \cite{PM0}. For a justification of some of the 
observations we refer to \cite{FZ}.
Following the path of Montgomery phenomenon (section \ref{sec:Montgomery}) we study the statistics of differences 
(deltas) $\gamma-\gamma'$ of imaginary parts of non-trivial Riemann zeros. But we don't normalize the differences, nor we restrict to 
nearby zeros as Montgomery did.
The elementary numerical study of the histogram of deltas $\gamma-\gamma'$, for $|\gamma|, |\gamma'|\leq T$ for $T$ large,  
shows some unexpected features (see Figure 3). 

We notice a uniform 
distribution on the first order, the Montgomery phenomenon located at $0$ at the next order, and at a lower order we have a deficit of deltas
at certain specific locations (see Figure 3). These locations are precisely the imaginary part of non-trivial zeros. One can say that 

\medskip
 
 \textbf{Riemann zeros repel their deltas.} 
 
\medskip

It is intriguing to note that the statistics of deltas is invariant by a translation, thus there is some global rigidity of the location of the 
non-trivial zeros. The statistical nature of the numerical computation shows that 

\medskip
 
 \textbf{Large Riemann zeros know about small Riemann zeros.} 

\medskip

\begin{figure}[h] 
\centering
\resizebox{9cm}{!}{\includegraphics{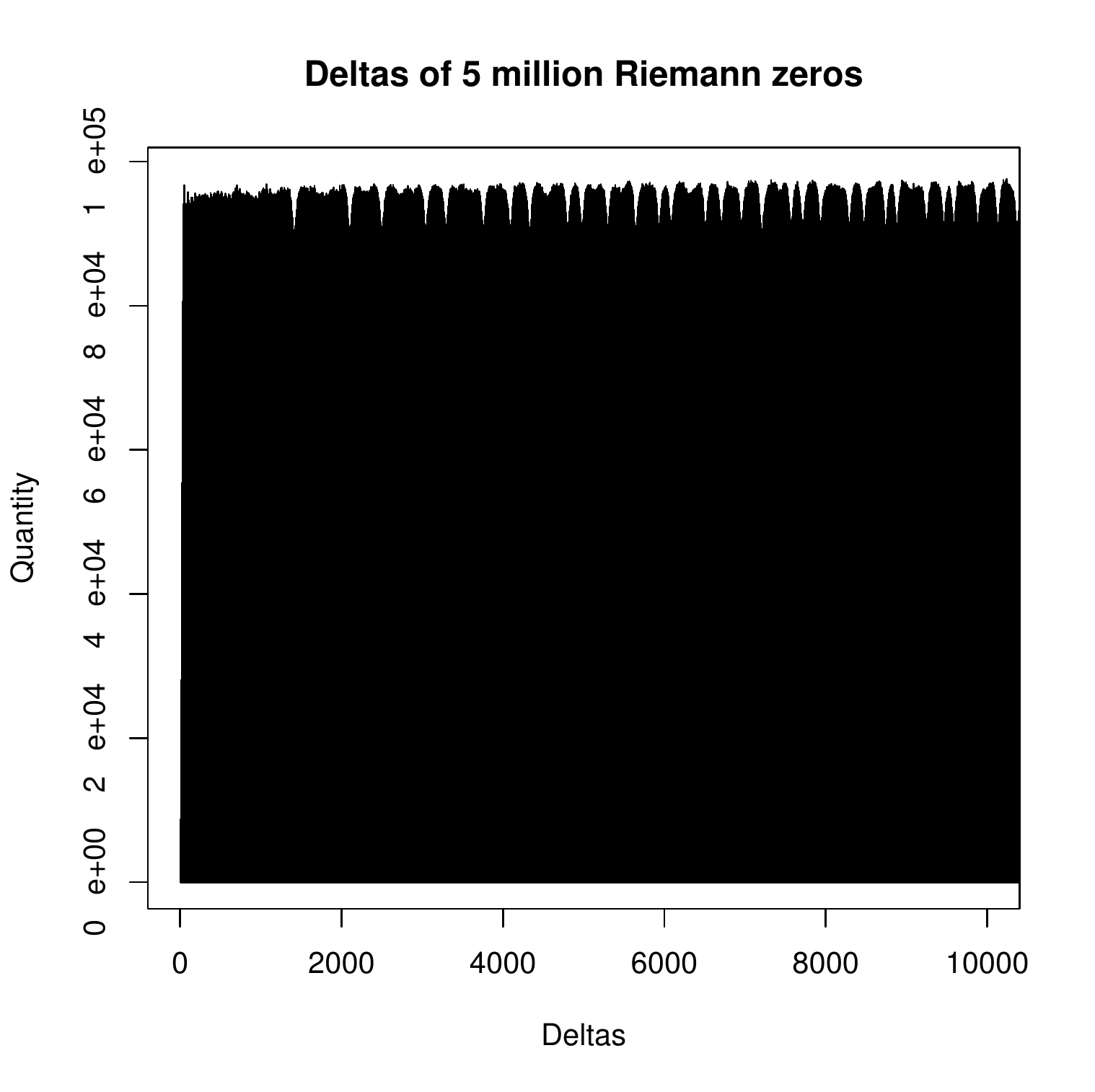}}    
\caption{Histogram of Riemann deltas.}
\end{figure}

 \medskip

 This type of phenomenon occurs also for other classes of zeta functions for which the Riemann Hypothesis is conjectured, for example for 
Dirichlet $L$-functions (see section \ref{sec:Dirichlet}). But in that case, in an even more surprising form since 
the statistics for deltas of non-trivial zeros show  
a deficit at the location of the Riemann zeros. We observe that 

\medskip

\textbf{Riemann zeros repel the deltas of $L$-functions zeros} 

\medskip

\begin{figure}[h] 
\centering
\resizebox{9cm}{!}{\includegraphics{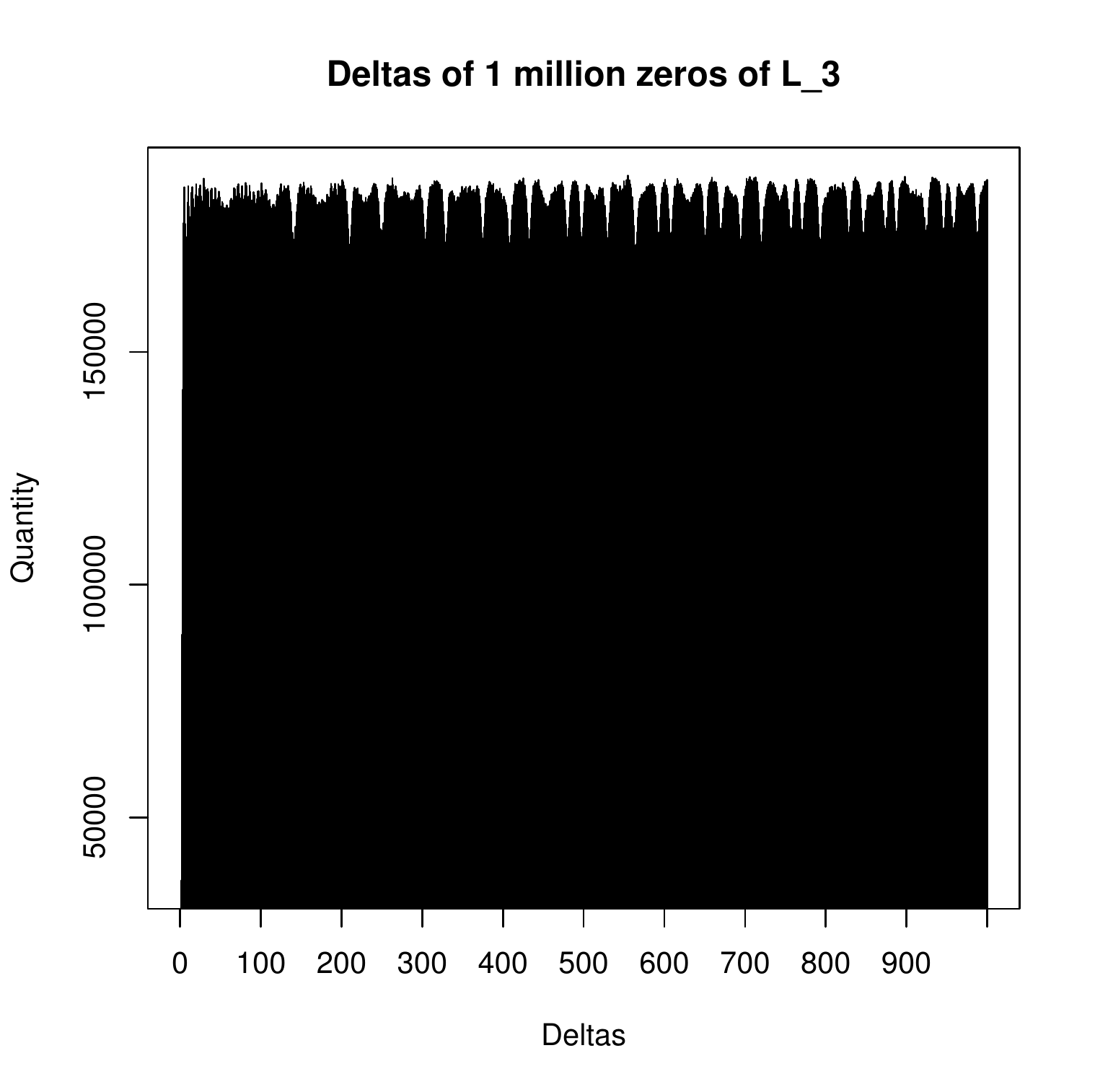}}    
\caption{Histogram of deltas of zeros of $L_3$.}
\end{figure}

 In Figure 4 the statistics are done with the zeros of the Dirichlet $L$-function (but there is nothing particular to this example)
$$
L_{\chi_3} (s)=\sum_{n=0}^{+\infty} \frac{\chi_3(n)}{n^s} \ ,
$$ 
where $\chi_3$ is the only character with conductor $3$, 
\begin{equation}
\chi_3(n)=
\begin{cases}
0, & \text{ if } n\equiv 0 \ \ \ [3] \\
1, & \text{ if } n\equiv 1 \ \ \ [3] \\
-1, & \text{ if } n\equiv -1 \ [3] 
\end{cases} 
\end{equation}
We also notice the Montgomery phenomenon.
These statistics allow to compute Riemann zeros using large zeros of Dirichlet $L$-functions. So we observe that

\medskip

\textbf{Non-trivial zeros of Dirichlet $L$-functions know about Riemann zeros} 

\medskip

There seems to be some confusion among physicists, and some believe that these observations 
are related to the pair correlation observations of the perturbative theory 
of the Quantum Chaos approach (compare for example the introductions of \cite{FZ} version 1 and version 3 from the ArXiv repository). 
But the true nature of these observations is completely different, and are more of an algebraic-arithmetic origin. 
The precise exact location of the zeros can be computed by performing statistics (and some more), which indicates 
that this is not related to any kind of perturbative theory (and indeed precise location of low zeros was never claimed by physicists). 
Also we can observe new phenomena that lacks of a physical interpretation. 

\medskip

More precisely, and surprisingly,  we can also mate zeros of two distinct Dirichlet $L$-functions: We do compute the statistics of the deltas $\gamma_1 - \gamma_2$ for 
$\gamma_1$, resp. $\gamma_2$, running over the non-trivial zeros of $L_{\chi_1}$, resp. $L_{\chi_2}$. We observe that the deficit of these deltas are located
at the imaginary parts of the non-trivial zeros of $L_{\chi_1 . \bar\chi_2}$. The intuition from Quantum Chaos is unable to explain this more general type 
of statistics since there isn't any type of ``pair correlation'' between energy levels of two different quantum systems.
It is also interesting to note that the Montgomery phenomenon, at the origin of quantum physical interpretations, does not arise when mating different $L$-functions
(different also from their conjugate).

\medskip

Moreover, we can also mate the zeros of $L_\chi$ with the poles of the Euler factors $\zeta_p (s)=(1-p^{-s})^{-1}$. Then we get as deficits 
the location of the imaginary parts of the poles of the local factors $L_{p,\chi}(s)=(1-\chi(p)p^{-s})^{-1}$.

\medskip

All of this reveals a hidden algebraic structure in the space of Dirichlet $L$-functions. The receipt to forecast all this results is the following operation, called the 
e\~ne-product, that can be defined formally at the level of Euler products. 
In the multiplicative group $\cA = 1+X\CC [X]$ we can define the e\~ne-product over $\CC$. Let $A,B\in \cA$, then if
\begin{align*}
 A(X) &=\prod_{\alpha} \left (1-\frac{X}{\alpha}\right ) \\
 B(X) &=\prod_{\beta} \left (1-\frac{X}{\beta}\right ) 
\end{align*}
we define
$$
A\star B (X) =\prod_{\alpha , \beta}  \left (1-\frac{X}{\alpha \beta}\right ) \ .
$$
Then it is easy to check that $(\cA, . , \star )$ has a commutative ring structure. Moreover, we can extend the e\~ne-product to $\cK$, the field of fractions of $\cA$ by 
$$
\frac{A}{B}\star \frac{C}{D}=\frac{(A\star C). (B\star D)}{(A\star D).(B\star C)}
$$
Observe that $A\star B = A^{-1}\star B^{-1}$ and
$$
(1-aX)\star (1-bX) =(1-aX)^{-1}\star (1-bX)^{-1}= 1-ab X \ .
$$
Now, given two Dirichlet series with an Euler product
\begin{align*}
 F(s) &=\prod_{p} F_p(p^{-(s-1/2)})  \\
 G(s) &=\prod_{p} G_p(p^{-(s-1/2)})  
\end{align*}
where the product is over all primes and $F_p, G_p \in \cK$, we can define their e\~ne-product $\bar \star$ as
$$
F\, \bar \star \, G (s)= \prod_p F_p \star G_p (p^{-(s-1/2)}) \ .
$$
The normalization is taken in that form because $z=e^{-s+1/2}$ is the correct variable. The space of Dirichlet series with the usual 
product and the e\~ne product then becomes a ring, where the multiplicative unit is the Riemann zeta function properly normalized. 
What these numerical observations show is that the e\~ne-product on arithmetic Dirichlet series has a deep interpretation 
at the level of the divisors (zeros and poles) of the functions. This divisor interpretation is clear for the simpler e\~ne-product of $\star$ in $\cK$.
But the interpretation for the divisors of Dirichlet functions is far more involved. The divisor of $F\, \bar \star \, G (s)$ is related with the 
additive convolution of the divisor of $F$ and the divisor of $G$ as for $\star$ in $\cK$, but this convolution does not give in general a discrete divisor. 
This is at the origin of the explanation of the statistics: The deltas of $L_{\chi_1}$ and $L_{\chi_2}$ point to the location of the zeros of $L_{\chi_1} \, \bar \star \, 
L_{\bar \chi_2}$. With this simple receipt we can forecast the numerical observations presented above. 

\medskip

The intuition on the $\bar \star$ product goes deeper. The fundamental ``unit equation'' ($1\times 1=1$)
$$
\zeta(s) \, \bar \star \, \zeta(s) =\zeta (s+1/2)^{-1} \ ,
$$
suggest that the Riemann hypothesis is true since $\zeta(s) \, \bar \star \, \zeta(s)$ will have no zeros in the half plane $\{\Re s > 1/2\}$. Indeed, a
non-trivial zero with $\Re \rho >1/2$ will induce the zero $ (\rho-1/2)+(\rho-1/2) -1/2 =2\rho-1/2$ (or $\rho+\bar \rho-1/2$) with real part $>1/2$. Note 
that this argument is independent of the normalization taken.
The situation for the e\~ne product in the space of arithmetic zeta functions is more far subtle, but we have here a good heuristic 
reason to believe in the Riemann Hypothesis.

\subsection{Adelic character.}

 Euler product formula and its relation to prime numbers shows that the Riemann zeta function has a deep 
 algebraic meaning. The reasons to believe in the Riemann Hypothesis are not only analytic. All the rich 
 generalizations of the conjecture to other transcendental zeta functions do arise for zeta functions with a 
 marked arithmetic character, more precisely with a well defined Euler product exhibiting local factors. This 
 is probably where the real difficulty lies: The truth of the conjecture must depend on a mixture of analytic 
 and arithmetic ingredients.

We have seen that the Euler product expansion naturally suggests the ``local" factors corresponding to each 
prime number,
\begin{equation*}
\zeta_p(s)=(1-p^{-s})^{-1} \ .
\end{equation*}
But the functional equation naturally introduces a more natural entire function $\xi$ which contains a new 
factor
\begin{equation*}
\frac{1}{2} s(1-s) \pi^{-s/2} \Gamma (s/2) \ .
\end{equation*}
How it should be interpreted? Tate in his thesis (1950, \cite{Ta} \cite{CF}) provides a uniform explanation 
to the new factor as the one corresponding to the infinite prime, and unifies the treatment of functional 
equations for various zeta and $L$ functions. Tate's approach is based on doing Fourier analysis at each local 
completion of $\QQ$ (or of a number field $K$, but we limit the exposition to $\QQ$ in this section). 
We recall Ostrowski's theorem \cite{Os} which classifies all absolute values over $\QQ$ up 
to topological equivalence.

\begin{theorem} \textbf{(A. Ostrowski)} Up to topological equivalence, all absolute values over the rational 
numbers are the ultrametric $p-$adic absolute value $|\cdot |_p$ and the archimedean absolute value 
$\abs{\cdot}_\infty.$
\end{theorem}

As usual, we note by $\QQ_p$ the completion of $\QQ$ for the $p$-adic absolute value $\abs{\cdot}_p$. The 
corresponding completion for $\abs{\cdot}_\infty$ is $\RR$. Tate does Fourier analysis on the id\`ele 
group
\begin{equation*}
 \AA^\times ={\prod_p}' \QQ_p^\times \ ,
\end{equation*}
where the restricted product means that only a finite number of coordinates are not in $\ZZ_p^\times$. We 
have an embedding $\QQ^\times \hookrightarrow \AA^\times$, and a decomposition
\begin{equation*}
 \AA^\times \approx \QQ^\times \times \RR_+^\times \times \hat \ZZ^\times \ ,
\end{equation*}
where 
\begin{equation}
 \hat \ZZ^\times =\varprojlim_n \left (\ZZ/n\ZZ\right )^\times \approx \prod_p \hat \ZZ_p \ .
\end{equation}
A quasi-character or Hecke character is a continuous complex character $\omega: \AA^\times \to \CC^\times$ which is trivial on $\QQ^\times$. These decompose into local characters $\omega=\otimes_v \omega_v$ and in each local 
completion $\QQ_v$ for $f_v \in \cS (\QQ_v)$  (the appropriate Schwartz class) we can define the Fourier-Mellin local transform 
\begin{equation*}
 Z_v(s, \omega_v, f_v)=\int_{\QQ_v^\times} f_v(x) \omega_v (x) |x|^s_v \ d_v^\times x \ ,
\end{equation*}
where $d_v^\times x$ is the multiplicative Haar measure normalize so that for $v_p$, $\ZZ_p^\times$ has volume $1$. We have also
a global Fourier-Mellin transform in the id\`eles for $f \in \cS (\AA)$
\begin{equation*}
 Z(s, \omega, f)=\int_{\AA^*} f(x) \omega (x) |x|^s \ d^\times x \ ,
\end{equation*}
and 
\begin{equation*}
 Z(s, \omega, f)=\prod_v Z_v(s, \omega_v, f_v) \ .
\end{equation*}
These transforms do have a meromorphic extension to all $s\in \CC$ and Fourier duality gives for them a 
functional equation for the substitution $s\mapsto 1-s$. From the complement factors we get the usual zeta
function and its local parts (see \cite{Ta} or \cite{Ku} for details). The explicit computation is
\begin{equation*}
\xi_p(s)=\int_{\QQ_p^\times}{\bf 1}_{\ZZ_p}(x)\abs{x}^s_p \ d^\times x=(1-p^{-s})^{-1}=\zeta_p(s) \ .
\end{equation*} 
where $d^\times x=\frac{1}{1-p^{-1}}\ \frac{dx}{\abs{x}_p}$. And for $p=\infty$, $d^\times x= \frac{dx}{\abs{x}}$ 
and ${\bf 1}_{\ZZ_p}$ is replaced by a Gaussian (the natural replacement would be by $\delta_0$, the Dirac at $0$, 
that is a limit of gaussians):
\begin{equation*}
\xi_\infty(s)=\int_\RR e^{-\pi x^2}\abs{x}^s \ \frac{dx}{\abs{x}}=\pi^{-\frac{s}{2}}\ \Gamma(\frac{s}{2}).
\end{equation*} 

The functions $\xi_p$ are the \textit{local factors} of $\zeta$. For $p=\infty$ is the real or 
infinite prime (this 
classical terminology should be replaced by ``the prime $p=1$'' as the Cramer function tells us 
clearly),  $\xi_\infty$ is called the \textit{local archimedean factor}. Then we reconstruct 
the global zeta-function as
\begin{equation*}
\xi(s)=\frac{1}{2}s(s-1)\xi_\infty(s)\prod_p\xi_p(s)
\end{equation*} 
The necessity of introducing a factor at infinity to obtain the functional equation is better understood in view of the analogy between number and function fields (see section \ref{sec:varieties}). As a matter of fact, $\mathrm{Spec}\ \ZZ$ is similar to the affine line $\mathbb{A}^1_k=\mathrm{Spec}\ k[T]$. The \textit{product formula} links together all different valuations
\begin{equation*}
\abs{x}_\infty\prod_p \abs{x}_p=1 \quad \forall x \in \QQ^\ast
\end{equation*} 
and this makes natural the simultaneous consideration on the same footing of all of them, and the hypothetical compactification $\overline{\mathrm{Spec}\ \ZZ}$. This is the point the view of Arakelov geometry. 

It remains to interpret the factor $s(s-1)$. An exotic interpretation was proposed by Y. Manin in \cite{Ma}: 
It would be the zeta-function of the projective line over the \textit{field with one element}.

\section{The universe of zeta-functions.} \label{sec:universe}

In order to better understand the true scope of the Riemann Hypothesis, 
we need to put it in context in the vast universe 
of zeta functions. These arise naturally in different fields of Mathematics, always with 
a rich associated structure.
Many times they encode deep arithmetical information. 
They provide unsuspected bridges and analogies between distinct fields that have proved to 
be the key points in the proof of central 
conjectures. The rich structure of such universe of zeta-functions is not well understood. 
Many mysterious relations and 
analogies appear with no easy explanation: 

\medskip

\textit{It's a whole beautiful subject and the Riemann zeta function is just the
first one of these, but it's just the tip of the iceberg. They are just the
most amazing objects, these L-functions -the fact that they exist, and
have these incredible properties are tied up with all these arithmetical
things- and it's just a beautiful subject. Discovering these things is like
discovering a gemstone or something. You are amazed that this thing
exists, has these properties and can do this. (J.B. Conrey)}

\medskip

It is very likely that a proof of the Riemann Hypothesis for 
Riemann zeta-function entails the unveiling of this hidden structure as a biproduct, and 
the proof is expected to work in a much larger class of zeta-functions. 

We start by the simplest and closer generalizations of Riemann zeta-function, and then we grow in 
generality and complexity. As we climb levels into increasing generality we notice a diminishing
knowledge of analytic properties.

\subsection{Dirichlet $L$-functions.}\label{sec:Dirichlet}
To any given Dirichlet character we can associate a Dirichlet $L$-function in the following way.
Let $\chi: \ZZ \longrightarrow \CC^*$ be a Dirichlet character defined modulo $q\geq 1$, 
that is, a $q$-periodic lift of a group morphism $(\ZZ/q\ZZ)^\times \longrightarrow \CC^*$ 
such that $\chi(n)=0$ when 
$(n,q)\not=1$. Thus we have the multiplicative property $\chi(nm)=\chi(n)\chi(m)$ and $|\chi(n)|=1$
when $\chi (n) \not= 0$. The character is \textit{primitive} when it is not induced by another 
character modulo $q'|q$. The 
\textit{conductor} is the minimal of such $q$'s. We define for $\Re s >1$,
\begin{equation*}
L(s,\chi)=\sum_{n=1}^{+\infty}\frac{\chi(n)}{n^s}=\prod_p (1-\chi(p)p^{-s})^{-1} \ .
\end{equation*} 
The multiplicative property of $\chi$ is what yields the Euler product. For $q=1$, $\chi=1$, we recover 
Riemann zeta-function, and for $q\geq 2$ and if the character is primitive and non-trivial ($\chi \not=1$), 
then $L(s, \chi )$ extends as an entire 
function over $\CC$, thus has no poles.

These  $L$-functions do satisfy a functional equation with symmetry respect to $\Re s=1/2$ 
when $\chi$ is a real valued character, or quadratic character, $\chi=\bar \chi$. Note that this case corresponds 
to have $L(s,\chi)$ real analytic. In general the 
functional equation only relates $L(s, \chi)$ and $L(s,\bar \chi)$.
If $\chi$ is a primitive character with conductor $q$ 
(we can reduce to this case modulo a finite number of Euler factors), then we 
have again a complement factor involving a $\Gamma$ factor satisfying the \textbf{DL} property. If
\begin{equation*}
\Lambda(s, \chi)= \pi^{-\frac{s+a_\chi}{2}} \Gamma\left (\frac{s+a_\chi }{2}\right ) L(s,\chi)
\end{equation*} 
where $a_\chi=0$ or $a_\chi=1$ depending whether $\chi(-1)=1$ (even character) or $\chi(-1)=-1$ (odd character). We denote by
$\tau(\chi)$ its Gauss sum,
\begin{equation*}
 \tau(\chi)=\sum_{x=0}^{q-1} \chi (x) e^{\frac{2\pi i}{q} . x} \ ,
\end{equation*}
and we define the \textit{root number} as $\epsilon(\chi)=q^{-1/2}\tau(\chi)$ for $\chi$ even, and $\epsilon(\chi)=-i q^{-1/2}\tau(\chi)$ for $\chi$ odd, then $|\epsilon (\chi)|=1$. Then we have
\begin{equation*}
\Lambda(\chi , s)=\epsilon(\chi) q^{1/2-s} \Lambda(1-s , \bar \chi ) \ .
\end{equation*}
When $\chi$ is a quadratic character, we have $\epsilon (\chi)=1$.

These $L$-functions where used by Dirichlet in  his proof that in any arithmetic sequence of integers 
$(a+nq)_{\n \in \ZZ}$, $(a,q)=1$, there are infinitely many primes, mimicking Euler's analytic proof of the existence of 
an infinite number of prime numbers. The result is equivalent to prove that when $\chi$ is non-trivial, then 
$L(1, \chi )\not=0$. Thus we find here another example where the location of zeros has a deep arithmetical meaning.

At this point, since for a general Dirichlet $L$-function we have no self-functional equation, 
we may wonder if the same 
intuition that has lead to formulate the Riemann hypothesis does makes sense for a non-quadratic 
character $L$-function. 
Observing that $Z(s, \chi)=L(s, \chi).L(s,\bar \chi)$ does indeed satisfy a 
self-functional equation, that it is 
a function in the 
\textbf{DL} class, with an Euler product, and that its divisor is just the 
union of both divisors, it is clear that again we can formulate the Riemann Hypothesis, that is 
indeed equivalent to the Riemann Hypothesis for both factors. 

\medskip

\textbf{Generalized Riemann Hypothesis.}  All non-real zeros of $L(s, \chi)$ lie on the 
vertical line $\Re s=1/2$.

\medskip

\subsection{Dedekind zeta-function.}\label{sec:Dedekind}

Dedekind's zeta-function is associated to a number field $K$ of degree $d=[K:\QQ]$.  
The ring of integers $\mathcal{O}_K$ is a Dedekind ring where we have unique factorization of integer ideals by prime 
ideals. Dedekind's zeta-functions is defined for $\Re s >1$ by 
\begin{equation*}
\zeta_K(s)=\sum_{\mathfrak{a}\subset \mathcal{O}_K} \frac{1}{N(\mathfrak{a})^s}=\prod_{\mathfrak{p}\subset 
\mathcal{O}_K}(1-N(\mathfrak{p})^{-s})^{-1} \ ,
\end{equation*} 
where $\mathfrak{a}$ runs over all ideals of the ring of integers, and $\mathfrak{p}$ over all prime ideals. 
For $K=\QQ$ we recover Riemann zeta-function $\zeta=\zeta_\QQ$.

Since every ideal in $\mathcal{O}_K$ has a unique factorization as a product of prime ideals, we have
\begin{equation*}
\zeta_K(s)=\prod_{\mathfrak{p}\subset \mathcal{O}_K}(1-N(\mathfrak{p})^{-s})^{-1}
\end{equation*} Hecke proved the meromorphic extension to $\CC$ with a simple pole at $s=1$, 
and the functional equation. If we write
\begin{equation*}
 \Lambda_K(s)= \pi^{r_1 s/2} (2\pi)^{r_2 s} \Gamma  (s/2)^{r_1} \Gamma (s/2)^{r_2} \zeta_K (s) \ ,
\end{equation*}
where $r_1$ and $r_2$ are the number of real and complex embeddings $K\hookrightarrow \CC$, then if
$D$ is the \textit{discriminant} of $K$,
\begin{equation*}
 \Lambda_K(s)= |D|^{1/2-s} \Lambda_K(1-s) \ .
\end{equation*}
Dedekind zeta-function $\zeta_K$ encodes in itself important arithmetic information. The best example is the 
analytic class number formula, that shows that we can read in the residue of $\zeta_K$ at $s=1$ important 
arithmetic information:
\begin{equation*}
 {\hbox{\rm{Res}}}_{s=1} \zeta_K(s)= \frac{2^{r_1} (2\pi)^{r_2} h(K) R}{w \sqrt{D}} \ ,
\end{equation*}
where $w$ is the number of roots of unity in $K$, $R$ is the regulator, and $h(K)$ the class number of $K$. 
The extension of this idea is what leads to Birch and Swinnerton-Dyer conjecture 
for $L$-functions of elliptic curves 
that is the subject of the 
exposition by V. Rotger in this same volume \cite{Ro}.

The Riemann Hypothesis is also conjectured in the same form for Dedekind zeta-functions.

\medskip

This construction was later generalized by Hecke (1917) with the introduction \`a la Dirichlet of characters in 
the ideal class group, and later (1918, 1920) with the introduction of ``Gr\"ossencharakter'' which is a 
character in the id\`ele class group. For these zeta and $L$-functions, Tate approach unifies the treatment of the 
local factors, meromorphic extension and functional equation.

\subsection{Artin $L$-functions.}

With Artin $L$-functions we jump another ladder in conceptual generality and we start loosing important analytic information. The reason is that Artin $L$-functions are defined through their Euler product, something that makes very difficult to prove meromorphic extension without an alternative analytic expression. 

Any Dirichlet primitive character of conductor $q$ is a character $\chi_{\Gal}$ on the group 
$\Gal (\QQ (e^{2\pi i/q})/\QQ) \approx (\ZZ/ q\ZZ )^\times$ (and hence on the profinite group 
$\Gal (\bar \QQ /\QQ)$), such that $\chi (p)=\chi_{\Gal}(\sigma_p)$, where $\sigma_p$ is the element of the Galois group given by the Frobenius automorphism $x\mapsto x^p$.  Moreover any $1$-dimensional character
on $\Gal (\bar \QQ /\QQ)$ is of this form. We can write the Dirichlet $L$-function as
\begin{equation}
 L(s, \chi)=\prod_p \left ( 1 - \chi_\Gal (\sigma_p) p^{-s}\right )^{-1} \ .
\end{equation}

Now, given a finite dimensional linear representation of an abstract group $G$, $\rho: G\to GL (V)$, its 
character is given by 
\begin{equation*}
 \chi(g) =\Tr \rho(g) \ .
\end{equation*}

When $K$ is a finite extension of a number field $k$ (of even a function field), and $\rho: G\to GL (V)$ 
a representation of the Galois group $G=\Gal(K/k)$ into the finite dimensional 
linear vector space $V$, and $\chi$ its 
character. For each prime ideal 
$\mathfrak{p}$ in $K$, the Frobenius automorphism $\sigma_{\mathfrak{p}}$ acts on $V^{I_\mathfrak{p}}$, where 
$I_\mathfrak{p}$ is the inertia group of $\mathfrak{p}$ in $G$ ($V^{I_\mathfrak{p}}=V$ when $\mathfrak{p}$ does 
not ramify). Then the local Artin $L$-function if defined by the characteristic polynomial of this action:
\begin{equation}
 L_p(s, \chi )=\Det \left (1-N(\mathfrak{p})^{-s} \sigma_{\mathfrak{p}} | V^{I_\mathfrak{p}} \right )^{-1} \ ,
\end{equation}
and then the global Artin $L$-function for finite places is given by the product of the local factors
\begin{equation}
 L (s, \chi )=\prod_p L_p(s, \chi ) \ .
\end{equation}
Note that the local factors here are no longer linear functions on $p^{-s}$, but are polynomial functions. In order to complete the picture, one introduces the archimedean factor that is a combination of 
\begin{align*}
 \Gamma_\RR (s) &= \pi^{-s/2} \Gamma (s/2) \\
 \Gamma_\CC (s) &= \Gamma_\RR (s) \Gamma_\RR (s+1)  \ ,
\end{align*}
in order to get $\Lambda(s, \chi )$.

With these definitions it is difficult to prove any meaningful analytic result, starting from the meromorphic extension. The only success is obtained when by algebraic means one could identify the $L$-function with 
others for which a more analytic (not from Euler's product) definition is at hand. The meromorphic 
extension, and the functional equation are indeed proved in two steps. First for the case of a one-dimensional 
representation: Class Field theory shows that $L(s, \chi)$ is a Hecke $L$-function, hence 
we get meromorphic extension, and the functional equation. Second, using a fundamental theorem of R. Bauer (1947)
on the classification of characters on a finite group, the general case is reduced to one-dimensional 
representations. We get an 
expression in terms of Dedekind zeta-functions. We get a functional 
equation of the form:
\begin{equation*}
 \Lambda (1-s, \chi)=W(\chi )  \Lambda (s, \chi) \ ,
\end{equation*}
where $W(\chi )\in \CC$, the \textit{Artin root number}, is of modulus $1$.
 
Again the Riemann Hypothesis is conjectured, and in this case follows from the Riemann Hypothesis 
for Dedekind zeta-functions. A. Weil generalized both Artin and Hecke zeta-functions by 
introducing a Gr\"ossencharakters in the construction.

\subsection{Zeta-functions of algebraic varieties.}\label{sec:varieties}

The well-known analogy between number fields and functions fields 
first started by R. Dedekind and H. Weber fundamental article \cite{DW}, has
always been a source of inspiration and testing. In general hard results for the number field case have been 
tested and proved  to be more accessible in the function field case. This also happens in the theory 
of zeta and $L$-functions. 

Along the lines of Dedekind's definition in the Number Field case, H. Kornblum used zeta-functions of 
algebraic varieties over a finite field to prove an analog of Dirichlet theorem on primes on 
arithmetic progressions. Later E. Artin continued this work and defined these zeta-functions in a 
more general context and studied the Riemann Hypothesis. Then Hasse and Schmidt 
proved the analogue of the Riemann Hypothesis for elliptic curves, and Weil finally proved it 
for general curves.
We define now zeta-functions for algebraic varieties over a number 
field $K$, or over a finite field $K=\FF_q$, $q=p^k$. 
We start by the finite field case $K=\FF_p$. The algebraic extensions of $\FF_p$ are 
the $\FF_{p^k}$ for $k=1,2,3,\ldots$ which are characterized dynamically as the set of periodic points of 
of period $k$ of the Frobenius automorphism $x\mapsto x^p$. Let $X_p$ be an algebraic variety defined over 
$\FF_p$ of dimension $n$. We define $N_k$ to be the number of points of $X_p$ in the field extension $\mathbb{F}_{p^k}$, and
\begin{equation*}
Z(t, X_p)=\exp \left (\sum_{k=1}^{+\infty}N_k \frac{t^{k}}{k} \right ) \ ,
\end{equation*}
or $\zeta (s, X_p)=Z(p^{-s}, X_p)$ in the more familiar $s\in \CC$ variable. Note that in this case, from the 
definition, the function $\zeta (s, X_p)$ is $i\omega_p^{-1}\ZZ$ periodic as the familiar local $\zeta_p$ factor. 
It is clear that this defines 
an analytic function in a right half plane. 

\medskip

\textbf{Example:}  For $X=\PP^1$, $Z(t, \PP^1)=(1-t)^{-1} (1-pt)^{-1}$.

\medskip

In 1949, inspired by his work on curves, Weil proposed a set of conjectures for these zeta-functions 
for $n$-dimensional algebraic varieties. The conjectures for curves were proposed in 1924 by Artin.

\begin{itemize}
 \item \textbf{Rationality.} The zeta-function $Z(t, X_p)$ is a rational function of $t\in \CC$, i.e. 
$\zeta (s, X_p)$ is a rational function of $p^{-s}$.
 \item \textbf{Functional equation.} We have $Z(\frac{1}{p^n t}, X_p)= \pm p^{nE/2} t^E Z(t, X_p)$, where $E$ is 
the self-intersection of the diagonal in $X_p \times X_p$.
 \item \textbf{Riemann Hypothesis.} We do have
\begin{equation*}
Z(t, X_p)=\frac{P_1(t)P_3(t)\ldots P_{2n-1}(t)}{P_0(t)P_2(t)\ldots P_{2n}(t)},
\end{equation*} where $P_0(t)=1-t, P_{2n}(t)=1-p^nt$ and for $1\leq j \leq 2n-1$, $P_j(t)=\prod_{k=1}^{B_j}(1-\alpha_{jk}t)$ where the $\alpha_{jk}$ are algebraic integers and $\abs{\alpha_{jk}}=p^{\frac{j}{2}}$. 
 \item \textbf{Betti numbers.} When $X_p$ is the reduction of a complete non-singular algebraic variety
in characteristic $0$, then the degree $B_j$ of $P_j$ is the $j$-th Betti number of $X$ as a complex manifold.
\end{itemize}

It is interesting to note that in this case, using the periodicity in $s\in \CC$, just the functional 
equation (second conjecture) plus the meromorphic extension to $\CC$ implies the rationality conjecture. Also 
the third conjecture is really a $\textbf{DL}$ property with $2n$ vertical lines, which shows again that the 
\textbf{DL} property arises naturally in this new context. 

In 1934 Hasse proved these results for elliptic curves and in 1940 Weil solved them for curves of 
arbitrary genus.
In 1960 B. Dwork, proved the rationality. Then M. Artin, A. Grothendieck, J.-L- Verdier and others 
developed $l$-adic \'etale cohomology theory.
This gave the polynomials $P_j$ \`a la Artin from the action of the Frobenius automorphism on cohomology. 
The functional equation then followed from Poincar\'e duality, and the degrees $B_j$ were identified 
as Betti numbers. 
Finally in 1974 P. Deligne proved the Riemann Hypothesis using the available technology of 
Lefschetz pencils and 
a positivity argument by Rankin in the theory of modular zeta-functions which is the key estimate, and that 
was obtained in this setting by a clever geometrical argument from an algebraic index theorem. As a corollary,
Deligne obtained a good estimate for the coefficients of modular functions associated to elliptic curves, 
in particular settling Ramanujan conjecture for the original modular function. The original goal of 
Rankin was precisely to obtain better estimates on the coefficients of modular functions. More precisely, 
Ramanujan $\tau$-function $\tau (n)$ is given by the modular form, which generates a modular zeta-function
(as explained in section \ref{sec:modular})
that is also the zeta function of an algebraic curve,
\begin{equation*}
\sum_{n\geq 1} \tau (n) q^n = q\prod_{n\geq 1} (1-q^n)^{24} = \eta (\tau )^{24}\ ,
\end{equation*}
where $\eta$ is Dedekind $\eta$-function. Ramanujan conjecture claims that for $p$ prime we have 
\begin{equation*}
\tau (p) \leq 2 p^{11/2} \ .
\end{equation*}

\bigskip

Now we discuss the case where $X$ is a non-singular algebraic variety defined over a number field, 
or just over $\QQ$ in order to simplify. 
For almost every prime, except for the 
so called ``bad reduction primes", the reduction $X_p$ of $X$ modulo $p$ is a non-singular algebraic 
variety $X_p$ over $\mathbb{F}_p$. We define the Hasse-Weil zeta function by the product (we avoid the finite
number of primes for which
$X_p$ may not be defined),
\begin{equation*}
Z_{X,\QQ}(s)={\prod}'_p Z_p(p^{-s}, X_p) \ .
\end{equation*} 
Note that \textbf{the zeta functions of $X$ over $\FF_p$ do appear as Euler local factors}. We should stress 
this point because we do have now a genuine transcendental function, and not a rational function as before. 
The mixture of rationally independent frequencies $(\omega_p)$ clearly shows that we are no longer facing a
vertically periodic function of the $s\in \CC$ variable.

\textit{Hasse conjecture} postulates the existence of a meromorphic extension. For the functional 
equation and the Riemann Hypothesis conjecture, one has to take into account the precise form of 
the Euler local factors (Serre, 1969), 
and usually is preferable to work separately with the zeta functions that are supposed to satisfy a minimal 
\textbf{DL} property (with only one vertical line for the location of the vertical divisor),
\begin{equation*}
 \zeta_j(s,X)={\prod}'_p P_j(p^{-s}, X_p)^{-1} \ .
\end{equation*}

\medskip
\textbf{Examples:}

\medskip

\textbf{(1)Projective space $\PP^1$.} For $X=\PP^1$, $\zeta_0(s, \PP^1)=\prod_p(1-p^{-s})^{-1}=\zeta(s)$ and we recover the usual Riemann zeta-function. 
Also $\zeta_2(s, \PP^1)=\prod_p(1-p p^{-s})^{-1}=\zeta(s-1)$, and 
\begin{equation*}
 Z_{X,\PP^1}(s)=\prod_p (1-p^{-s})^{-1} (1-p p^{-s})^{-1}=\zeta(s) \zeta(s-1) \ .
\end{equation*}

\medskip

\textbf{(2) Elliptic curves.} We refer to the text of V. Rotger in this same volume \cite{Ro} for examples 
and a detailed discussion. We just mention that the key point in the resolution of Fermat problem is the proof 
of the Taniyama-Weil conjecture. This conjecture postulates the modularity of the zeta-function of elliptic 
curves, i.e. that these $L$-function come also from modular forms in the sense of section \ref{sec:modular}, 
and this yields many previously unknown analytic properties. This is a very good example of the 
kind of bridges between distinct parts of the zeta-function universe  
that provides an insight in many problems.

\bigskip

\textbf{Comments.}

At this point it should be clear that the zeta-functions of varieties over finite fields 
represent nothing more than local Euler factors of the Hasse-Weil transcendental zeta-functions. Despite the 
difficulties and technicalities involved in the proof of the Riemann Hypothesis for these local 
zeta functions, coming back to the discussion in the first section, the result is equivalent in this 
context to Euler's proof that the sine function has no complex zeros, i.e. that the local factor 
$\zeta_p(s)$ satisfies the Riemann Hypothesis. How can this result represent a significant advance 
towards the Riemann Hypothesis? In no way in the authors opinion. The knowledge of the Riemann Hypothesis 
for $\zeta_p(s)$ hasn't provide yet any clue about how to go about the transcendental aspects of the Riemann Hypothesis for the 
transcendental zeta function $\zeta(s)$. One may agree with E. Bombieri's claim \cite{Bo} that 

\medskip

\textit{``Deligne's theorem surely ranks as one of the crowning achievements of 20th century mathematics.`" (E. Bombieri)}

\medskip

because of the important techniques and inspiring arguments used in the proof, but there is no reason to claim 
that this represents a significant advance towards the Riemann Hypothesis for transcendental 
zeta-functions. The Riemann Hypothesis for rational zeta functions (rational in the 
appropriate variable, or periodic 
in $s$), are ``toy zeta-functions'' where one can test the Conjecture, and the algebraic arguments, but they 
cannot give any clue about the intrinsic transcendental difficulty attached to 
the problem. It is natural that pure algebraists may downsize or overlook these 
difficulties. But one should be very careful 
in how things are presented and what claims are made in order to get the correct 
intuition on the real problem.

It is interesting to note that another example of ``toy zeta-functions'' is 
provided by dynamical zeta-functions associated to hyperbolic 
dynamical systems (see section \ref{sec:dynamical}). 
As we will see in section \ref{sec:dynamical}, it is possible to encode the 
dynamics with a finite Markov partition, which provides a combinatorial encoding of periodic orbits, and 
a combinatorial framework that is used to prove the rationality of the zeta-function. This result 
is of the same nature than Deligne's result.

\subsection{$L$-functions of automorphic cuspidal representations.}

Inspired by the adelic ideas of Tate, in the seventies there was a 
development for a general theory of a very large class 
of $L$-functions associated to automorphic cuspidal representations of 
$GL(n, \AA )$ where $\AA$ is the ad\`ele
group of $\QQ$, that are supposed to generalize all previous zeta and $L$-functions. Such a representation $\pi$ is equivalent to $\otimes_v \pi_v$ where $\pi_v$ is an irreducible unitary of $GL (n , \QQ_v))$. We 
will give an sketchy presentation because of lack of space in order to introduce all the terminology.

When $v_p$ is a finite prime, one determines from $\pi_p$ a local factor $L(s, \pi_p)$ which is a polynomial
on $p^{-s}$ of degree $n$. Then the associated $L$-function is
\begin{equation*}
L (s, \pi)=\prod_{p} L(s, \pi_p) \ .
\end{equation*}
For $v=\infty$, then $\pi_\infty$ determines $n$ parameters $(\mu_{j,\pi_\infty})_{1\leq j \leq n}$ and
\begin{equation*}
L (s, \pi_\infty)=\prod_{j=1}^{n} \Gamma (s-\mu_{j,\pi_\infty} ) \ .
\end{equation*}
Then the function
\begin{equation*}
\Lambda (s, \pi)= L (s, \pi_\infty) \ . \ L (s, \pi) \ ,
\end{equation*}
satisfies the functional equation
\begin{equation*}
\Lambda (s, \pi)= \epsilon_\pi N_\pi^{1/2-s}  \Lambda (s, \tilde \pi)\ ,
\end{equation*}
where $|\epsilon_\pi |=1$, $N_\pi \geq 1$ is an integer, and $\tilde \pi$ is the contragredient representation.

A general conjecture of Langlands is that these $L$-functions generate multiplicatively all previous $L$-functions (Dirichlet, Dedeking, Hecke, Artin, Hasse-Weil,...). And we expect the Riemann Hypothesis to hold in this 
class of $L$-functions:

\medskip

\textbf{Grand Riemann Hypothesis.} The zeros of $\Lambda (s, \pi )$ all have real part $1/2$.

\medskip

Again, as for Artin $L$-functions, one of the weak points of the general 
definition of $L$-functions associated to 
automorphic representations is that the definition by the Euler product is analytically unusable. It is 
necessary already from the very definition to prove that the $L$-function has a meromorphic extension to
$\CC$.

\subsection{Classical modular $L$-functions.}\label{sec:modular}

We recall that Riemann gave an integral formula expressing Riemann $\zeta$-function from the 
classical $\theta$-function, and the functional equation for $\zeta$ resulted from the modular functional 
equation for $\theta$. This procedure generalizes to give a large family of zeta-functions associated to 
certain modular forms. These form the class of modular $L$-functions that have good analytic properties. This 
time we are taking a different path and are not escalating in generality.

Let $\mathbb{H}=\{z \in \CC: \Im z>0\}$ be the upper half plane and let 
$H$ be a subgroup of $\mathrm{SL}(2,\ZZ)$ containing a congruence subgroup 
\begin{equation*}
\Gamma(N)=\left\{ \left( \begin{array}{cc}
a & b  \\
c & d \end{array} \right) \equiv \left( \begin{array}{cc}
1 & 0  \\
0 & 1 \end{array} \right) [N] \right\} \ .
\end{equation*}

\begin{definition} A \textit{meromorphic modular form} of weight $k$ relative to $H$ is a meromorphic 
function $f: \mathbb{H} \longrightarrow \overline{\CC}$ such that
\begin{equation*}
f\left (\frac{az+b}{cz+d} \right )=(cz+d)^k f(z) \ \text{ for all matrices } \ \left( \begin{array}{cc}
a & b  \\
c & d \end{array} \right) \ \text{ in } \ H \ ,
\end{equation*}
and $f$ is meromorphic at the cusps, that is, when H acts on $f$, the resulting function has a 
meromorphic expansion in $q^{1/N}$ at $q=0$ for the variable $q=e^{2 \pi i z}$, that is
\begin{equation*}
 (cz+d)^k f\left (\frac{az+b}{cz+d} \right )=\sum_{n=-r}^{+\infty } a_n \ q^{n/N} \ .
\end{equation*}

We say that $f$ is a \textit{modular form}, resp. \textit{parabolic modular form}, 
when it is holomorphic at the cusps (always $r\geq 0$), resp. vanishes at the cusps ($r\geq 1$).
\end{definition} 

In order to simplify, we limit the discussion to $N=1$.
Let $M_k$, resp. $S_k$, be the vector spaces of modular, resp. modular parabolic, forms 
for $H= \mathrm{PSL}(2,\ZZ)$. 
The space $M_k$ is finite dimensional and its dimension can be computed explicitly. We have
$\dim M_k=[k/12]$ if $k\equiv 2 [12]$, and $\dim M_k=[k/12]+1$ otherwise. Moreover we have, 
\begin{equation*}
M_k=S_k\oplus \CC G_k \ ,
\end{equation*} 
where $G_k$ is the Eisenstein series. 

From $f\in S_k$ and more precisely, its expansion at infinite, we construct an $L$-function by
\begin{equation*}
L(f,s)=\sum_{n=1}^{+\infty}\frac{a_n}{n^s} \ .
\end{equation*}
In order to have good properties for this $L$-functions (Euler product, functional equation,...) 
we need to restrict our class of modular forms to those with good arithmetic properties. The selection 
is done taking eigenfunctions of the Hecke operators. 

For $n\geq 1$, we define the Hecke operator $T_n: M_k \longrightarrow M_k$ or $T_n: S_k \longrightarrow S_k$ by

\begin{equation*}
T_n f(\tau)=n^{k-1}\sum  (c\tau+d)^{-k} f \left (\frac{a\tau+b}{c\tau+d} \right ),
\end{equation*} 
where the sum runs over the matrices 
$\left (\begin{array}{cc}
a & b  \\
c & d \end{array} \right )\in \mathrm{SL}(2,\ZZ )\\ \mathcal{M}_n$
where $\mathcal{M}_n$ is the set of integral matrices of order two having determinant $n$. 

Since
\begin{equation*}
T_nT_m=\sum_{d|(n,m)}d^{k-1}T_{\frac{nm}{d^2}}  \ ,
\end{equation*} 
when $n$ and $m$ are coprimes, $T_nT_m=T_mT_n$. Then $T_n$ and $T_m$ diagonalize in the same bases. A Hecke 
form is an eigenvector $f\in S_k$ for all $T_n$. If $f(z)=\sum_{n=1}^{+\infty}a_nq^n$, then $a_1\neq 0$, so 
it can be normalized to $a_1=1$. The coefficients are algebraic numbers $a_n \in \overline{\QQ}$ with degree bounded 
by the dimension of $S_k$. Moreover we have the multiplicative relations for the $a_n$:
\begin{equation*}
a_n a_m=\sum_{d|(n,m)}d^{k-1}a_{\frac{nm}{d^2}} \ .
\end{equation*} 
One can find a bases of Hecke forms for $M_k$. Now if $f$ is a Hecke form, the associated $L$-function has 
an Euler product that comes from the multiplicative relations for the $(a_n)$:
\begin{equation*}
L(s, f)=\sum_{n=1}^{+\infty}\frac{a_n}{n^s}=\prod_p(1-a_pp^{-s}+p^{k-1}p^{-2s})^{-1} \ .
\end{equation*} 
Then $L(s,f)$ converges for $\Re s>k$, and for $\Re s > \frac{k}{2}+1$ when $f$ is a cuspidal form. 
One can prove the meromorphic extension to all of $\CC$. It is an holomorphic function for $f$ cuspidal and has a simple pole at $s=k$ otherwise. We have the functional equation:
\begin{equation*}
(2\pi)^{-s}\Gamma(s)L(s,f)=(-1)^{\frac{k}{2}}(2\pi)^{s-k}\Gamma(k-s)L(k-s, f) \ .
\end{equation*} 

\subsection{$L$-functions for $GL(n, \RR)$.}

A vast generalization of the precedent modular $L$-functions was started by H. Maass in 1949 
when he noticed that we 
don't really need classical modular forms in order to generate $L$-functions, but non-holomorphic automorphic forms are enough.

More precisely, one looks for eigenfunctions of the Laplace operator for the Poincar\'e metric in the upper 
half plane which is invariant by the action of $SL(2, \ZZ)$:
\begin{equation*}
\Delta = -y^2 \left (\frac{\partial^2}{\partial x^2} + \frac{\partial^2}{\partial y^2} \right ) \ .
\end{equation*}
The function $I_s(z)=y^s$  is an eigenfunction with eigenvalue $s(1-s)$, and averaging over the action 
of $SL(2,\ZZ )$ we do get the Eisenstein series:
\begin{equation*}
E(z,s) = \frac{1}{2}\sum_{c,d\in\ZZ ; (c,d)=1} \frac{y^s}{|cz+d|^s} \ .
\end{equation*}
Its Fourier expansion can be computed:
\begin{equation*}
E(z,s) = y^s+\varphi (s) y^{1-s} + \frac{2 \pi^s \sqrt{y}}{\Gamma (s) \zeta (2s)} \sum_{n\in \ZZ^*} 
\sigma_{1-2s}(n) |n|^{s-1/2} \ K_{s-1/2} (2\pi |n| y) \ e^{2 \pi i n x}\ ,
\end{equation*}
where 

\begin{equation*}
\varphi(s)=\sqrt{\pi } \frac{\Gamma(s-1/2)}{\Gamma(s)} \frac{\zeta(2s-1)}{\zeta(2s)} \ ,
\end{equation*}

\begin{equation*}
\sigma_s(n)=\sum_{d | n} d^s \ ,
\end{equation*}
and $K_s$ is the Bessel function
\begin{equation*}
K_s(y)=\frac{1}{2} \int_0^{+\infty} e^{\frac{1}{2} y (u+1/u)} \ u^s \ \frac{du}{u} \ .
\end{equation*}
Thus we see that from these Fourier coefficients we can construct the $L$-function associated to 
the Eisenstein series
\begin{equation*}
L_{E( . ,w)}(s)=\sum_{n\geq 1} 
\sigma_{1-2w} (n) n^{w-1/2-s} =\zeta(s+w-1/2) . \zeta(s-w+1/2) \ .
\end{equation*}

 This motivates to look for a general definition for Maass forms:
 
\begin{definition} Let $L^2 (\HH )$ be the space of square integrable functions with respect to the Poincar\'e
volume.
A Maass form is $f\in L^2 (\HH )$ such that
\begin{itemize}
\item $\forall \gamma \in SL(2,\ZZ )$, we have $f(\gamma . z) =f(z)$.
\item $\Delta f = \nu (1-\nu ) f$.
\item $\int_0^1 f(z) \ dz=0$.
\end{itemize}
\end{definition}

If $f$ is a mass form, then it has a Fourier expansion of the form
\begin{equation*}
  f(z)=\sum_{n\in \ZZ^*} a(n) \ \sqrt{2 \pi y} \ K_{\nu-1/2} (2\pi |n| y) \  e^{2\pi i n x} \ ,
\end{equation*}  
  
We can define commuting Hecke operators $(T_n)_{n\geq 1}$ commuting with the Laplacian and acting on the space of mass forms. If $f$ is a mass form 
normalized with $a(1)=1$ and being the eigenvalue of all Hecke operators, then for $n=1,2,3,\ldots$
\begin{equation*}
T_n(f)= a(n) f \ ,
\end{equation*}
and the $(a_n)$ do satisfy multiplicative relations. If $f$ is odd or even and 
we define the associated $L$-function by
\begin{equation*}
L_f(s)=\sum_{n=1}^{+\infty} a(n) \ n^{-s} \ ,
\end{equation*}
the multiplicative relations give the Euler product:
\begin{equation*}
L_f(s)=\prod_p (1-a(p) p^{-s} + p^{-2s} )^{-1}\ .
\end{equation*}
One can show that these $L$-functions have a holomorphic continuation to all of $\CC$ and 
satisfy the following functional equation:
\begin{equation*}
\Lambda_f (s)=(-1)^{\epsilon} \Lambda_f (1-s) \ ,
\end{equation*}
where
\begin{equation*}
\Lambda_f (s)=\pi^{-s} \Gamma\left ( \frac{s+\epsilon -1/2+\nu}{2} \right ) +\Gamma\left ( \frac{s+\epsilon +1/2-\nu}{2} \right ) \ L_f(s)\ ,
\end{equation*}
where $\epsilon=0$ if $f$ is even, and $\epsilon =1$ if $f$ is odd.

Through the Iwasawa decomposition in $GL(n,\RR)$ we can define a generalization of the upper half plane 
$\HH$,
\begin{equation*}
\HH_n = GL(n, \RR)/(O(n, \RR). \RR^*) \ ,
\end{equation*}
and we can consider the algebra of $SL(n, \ZZ)$ invariant differential operators, Hecke operators 
and the corresponding Maass for 
$SL(n, \ZZ)$. Then through their Fourier coefficients, R. Godement and H. Jacquet defined their associated $L$-functions that have an Euler product whose local factors are polynomials of degree $n$ in $p^{-s}$. 
The Rankin-Selberg convolution 
operator can be defined, and has better analytic properties in the $GL(2, \RR)$ theory that in the higher 
dimensional $GL(n, \RR )$ theory where it is more formally defined.

This vast generalization lead to R.P. Langlands to a set of conjectures known as ``The Langlands Program".

\subsection{Dynamical zeta-functions.}\label{sec:dynamical}

Zeta-functions encode in a wonderful way arithmetical information. They can also encode dynamical and geometrical 
information as we see now. 

Following the path of the definitions of zeta-functions over algebraic varieties, where we noted the role played by the dynamics of the Frobenius map, we can define in a similar way the zeta-function of a dynamical system $(f, X)$,  $f: X\to X$. 
For $j=1,2,3,\ldots$, let $N_j(f)$ be the number of fixed points of the $j$-th iterate $f^j$, that is the number of 
distinct solutions to $f^j(x) =x$ (we will assume that this number is finite). Then the zeta-function associated to $f$ is
\begin{equation*}
 Z_f(t)=\exp \left (\sum_{j=1}^{+\infty}  \frac{N_j(f)}{j} \ t^j \right ) \ .
\end{equation*} 
If the number of fixed points grow at most in an exponential way, this defines a holomorphic function in a neighborhood 
of $t=0$. 

\medskip

\textbf{Fundamental example: Subshift of finite type.}

We consider the dynamics of a shift of finite type. We denote $\Sigma(l)= \{1,2,\ldots , l\}^\ZZ$ the compact 
space endowed with the product topology. An element $x \in \Sigma(l)$ is a sequence of the $l$ symbols of the 
alphabet, $x=(x_j)_{j\in \ZZ}$. The space is metrisable and a Cantor set. The shift homeomorphism  is the map
$\sigma : \Sigma (l)\to \Sigma (l)$ defined by 
\begin{equation*}
 \sigma (x)_j =x_{j+1} \ .
\end{equation*}

Let $A=[a_{ij}] \in M_{l\times l} (\{0,1\})$ be a square $l\times l$ matrix of zeros and ones. We can associated to 
$A$ the compact subspace
\begin{equation*}
 \Sigma_A=\{ x\in \Sigma (l) ;  \forall j\in \ZZ, \ a_{x_j x_{j+1}} =1 \} \ .
\end{equation*}
The shift $\sigma$ leaves $\Sigma_A$ invariant. 
The subshift of finite type associated to $A$ is the is the dynamical system $(\sigma_A , \Sigma_A)$ with
$\sigma_A : \Sigma_A \to \Sigma_A$ where $\sigma_A =\sigma_{/\Sigma_A}$ is the restriction of the shift.

This type of systems arises continually in dynamics when we want to code a dynamical system. If $(X_j)_{1\leq j\leq l}$ is a partition of $X$ (Markov partition), and $f:X\to X$ is an invertible 
dynamical system on $X$, we can associate a transition matrix $A$ by
defining $a_{ij}=1$ if $f(X_i)\cap X_j \not=\emptyset$ or $a_{ij}=0$ otherwise. Then we can encode the orbit of any $x\in X$ by associating to $x$ the sequence of symbols $\varphi(x)$ with 
\begin{equation*}
 f^j(x) \in X_{\varphi(x)_j} \ .
\end{equation*}
When the dynamical system is not invertible, we can consider only the lateral shift, or we can also 
construct an invertible dynamical system taking an inverse limit by a classical construction.

We can combinatorially enumerate the number of periodic orbits of the subshift of finite type in a very simple way:

\begin{lemma}
 For $n\geq 0$, we have
\begin{equation*}
N_n(\sigma_A )=\Tr A^n \ .
\end{equation*}
\end{lemma}

From this we get that the associated zeta-function is rational:
\begin{equation*}
 Z_{\sigma_A} (t)=\exp \left (\sum_{j=1}^{+\infty}  \frac{\Tr A^j}{j} \ t^j \right )= \left (\Det (I-t A) \right )^{-1} \ . 
\end{equation*}

There are some systems that do have a good coding and for which we can prove the rationality 
of the zeta-function 
in full generality. These are the hyperbolic dynamical systems (see \cite{Sh}). For example, the linear
map given by the matrix $\left (\begin{array}{cc}
2 & 1  \\
1 & 1 \end{array} \right )$ on the two dimensional torus $\TT^2$ is hyperbolic. Such systems are stable 
under perturbations. A more general notion is of Axiom A difeomorphism 
$f:M\to M$ of a compact manifold $M$, which
means the the non-wandering set $\Omega(f)$ of $f$ is hyperbolic and periodic orbits are dense on it. Then we 
can construct a Markov partition of $\Omega(f)$ with rectangular regions delimited by pieces of stable and 
unstable manifolds, using the existing local product structure: For two nearby point $x,y \in \Omega(f)$, we can 
grow a small local stable manifold from $x$, $W_\epsilon^S(x)$, and a small  local
unstable manifold from $y$, $W_\epsilon^U(y)$, such that they intersect in a single point 
$[x,y]=W_\epsilon^S(x) \cap W_\epsilon^U(y)\in \Omega (f)$. Except for the 
orbits hitting the common boundaries, the orbit 
is uniquely determined by its coding. Thus, in such a way we can count the periodic orbits, except for the 
redundancy on the orbits going through boundaries of the Markov partition. This allows to prove, by pure 
combinatorial enumeration, that the zeta function is rational.

It may seem surprising that for example an hyperbolic diffeomorphism of the torus $\TT^2$ 
has such a rational function, considering that rationality should come 
from an algebraic phenomena. The mystery disappears when one learns 
that such an arbitrary diffeomorphism is always topologically equivalent to a linear map, and the count of
periodic orbits is clearly the same for both systems.

\bigskip

We have discussed dynamical zeta-functions for discrete iteration. We can also define zeta-functions for 
continuous dynamical systems (flows) by measuring (with a natural or auxiliary 
Riemannian metric) the length of periodic 
orbits. For example, Selberg zeta-function for the geodesic flow of a compact surface of constant negative
curvature is of this nature. Selberg zeta function encodes the length spectrum, i.e. the set of length of 
primitive simple geodesics.


\subsection{The Selberg class.}

In 1992 A. Selberg proposed the more general class of zeta-functions for which he expected the Riemann Hypothesis to hold.

The Selberg class $\cS$ is composed by those Dirichlet series \footnote{There is a much larger class of functions that are not Dirichlet series for which to expect the Riemann hypothesis.}  
\begin{equation*}
F(s)=\sum_{n=1}^{+\infty} a_n n^{-s}\ ,
\end{equation*}
that are absolutely convergent for $\Re s>1$ that satisfy the following four conditions:
\begin{itemize}
\item \textbf{Analyticity.} There is $m\geq 0$ such that $(s-1)^m F(s)$ is an entire function.
\item \textbf{Ramanujan Conjecture.} $a_1=1$ and $a_n = o(n^\epsilon )$ for any $\epsilon > 0$.
\item \textbf{Functional equation.} There is a factor of the form
\begin{equation*}
\gamma_F (s) = \epsilon Q^s \prod_{j=1}^k \Gamma (w_j s+\mu_j) \ , 
\end{equation*}
where $|\epsilon |=1$, $Q > 0$, $w_j > 0$, $\Re \mu_j \geq 0$ such that
\begin{equation*}
\lambda (s)= \gamma_F(s) F(s) \ ,
\end{equation*}
satisfies the functional equation
\begin{equation}
\Lambda (s)=\overline{\Lambda (1-\bar s )} \ .
\end{equation}
\item \textbf{Euler Product.} We have for $\Re s > 1$
\begin{equation*}
F(s) = \prod_{p}F_p(s) \ ,
\end{equation*}
and
\begin{equation*}
\log F_p(s) = \sum_{k=1}^{+\infty } b_k p^{-ks} \ ,
\end{equation*}
and for some $\theta \in ]0,1/2[$, we have $b_{p^k} = \cO (p^{k\theta})$.
\end{itemize}

On top of the Riemann Hypothesis, Selberg made a set of conjectures for this class of Dirichlet 
series known by the name of \textit{Selberg Conjectures.}

\subsection{Kubota-Leopold zeta-function.}\label{sec:Kuwota-Leopold}
These zeta and $L$-functions are of a very different nature than the previous examples. This time, the generalization goes into a different direction: These are zeta-functions of a $p$-adic variable.

The starting point is the observation of the remarkable rational values given by Bernoulli numbers, taken by $\zeta(s)$ at the negative integers, for $n\geq 1$,
\begin{equation*}
\zeta(1-n)=-\frac{B_n}{n} \in \QQ \ .
\end{equation*} 
Let $p\neq 2$ (for $p=2$ the formulas change) be  a prime number and fix $u\in \ZZ/(p-1)\ZZ$. 
When $n_k \to s \in \ZZ / (p-1) \ZZ \times \ZZ_p$ as $k \to +\infty$, with the restriction $n_k \equiv u$ mod $p-1$, one can prove 
(using Kummer congruences for Bernoulli numbers) that the sequence $(\zeta(1-n_k))_k$ is converging 
in the $p$-adic field $\QQ_p$. Indeed there exists an analytic function in $\QQ_p$ defined by
\begin{equation*}
\zeta_{p,u}(s)=\lim_{k\to +\infty} \zeta (1-n_k) (1-p^{n_k}) \ ,
\end{equation*} 
which interpolates the values taken by $\zeta(s)(1-p^s)$ at the $s=1-n_k$. For $u$ even we get $\zeta_{p,u}=0$.

We have the formulas, for odd $u$ and $n\leq 0$, $n\equiv u [p-1]$
\begin{equation*}
 \zeta_{p,u}(n)=(1-p^{-n}) \zeta (n) \ , 
\end{equation*}
and for $n\geq 1$, $n\equiv u [p-1]$,
\begin{equation*}
 \zeta_{p,u}(n)=\frac{\Gamma(n)}{(2\pi i)^n}(1-p^{n-1}) \zeta (n) \ . 
\end{equation*}
We should note that $u+(p-1) \NN$ and $u-(p-1) \NN$ are dense in $\ZZ_p$.

\section{Conclusion and comments on the Riemann Hypothesis.}

\subsection{Is it true?}

The first question one has to ask about an open conjecture is if we do actually believe it. 
In the case of the Riemann Hypothesis there seems to be an overwhelming opinion in favour of it. And this, 
despite that the vast majority of mathematicians that have an opinion don't have a clue on how to go about it. 
Also despite some notable exceptions of unbelievers. Indeed, it appears that some of the best specialists that 
have spend considerable effort into it, at the end of their life start to become sceptic. We can just mention 
the best known of them, J.E. Littlewood. At first this seems quite troubling. But knowing human nature, 
probably we can only take this as another indication that the Riemann Hypothesis is certainly true. One can 
notice that most (all?) of the attempts to resolve the questions have tried to prove and not disprove it. The 
only attempts to disprove it seem to have been numerical. We will not even discuss the possibility of the 
Riemann Hypothesis to be non-decidable: It is obviously a genuine well posed beautiful problem. 

The arguments in favour of the Riemann Hypothesis that one finds in the literature are weak.

\begin{itemize}
 \item \textbf{Numerical evidence.} As some authors have observed, the range at the scope of present 
 computer power is probably not  large enough to draw any conclusion in view that ranges 
 where $\log \log T$ is large should be 
 inspected (see \cite{Iv2}).
 
 \item \textbf{Proportion of zeros in the critical line.} Goes into the right direction, but says nothing about the totality of zeros.
 
 \item \textbf{Weak density of zeros away from the critical line.} Same comment as before.
 
 \item \textbf{True for toy models.} As the zeta function of algebraic varieties over finite fields, or hyperbolic Dynamical Systems. Obviously these go into the right direction, but miss the central analytic difficulty for transcendental functions. This explain why these results are tractable by purely geometrico-algebraic, or combinatorial methods, but the transcendental problem seems totally out of reach.
 
 \item \textbf{Spectral Hilbert-Polya approach.} It is certainly a very appealing proposition, but does not give direct support unless we have a proposition for the unitary operator that will have the non-trivial zeros as spectrum. The dream of finding a physical hamiltonian-like system with this spectrum of energies is very unlikely to succeed. In view of the universe of zeta-functions, the system for Riemann zeta-function would be a very basic and almost canonical one. Thus easy to find. So far, no such system has been proposed. It is very likely that it does not exist. Most probably the part of the hidden structure of the universe of zeta-functions is shared by mechanical systems, which would explain the observed analogies.
 
\end{itemize}

\subsection{Why is it hard?}

As most of the thought problems, it is hard because it contains in a non-trivial way a 
mixture of analytic and arithmetic elements. One can probably argue the same way 
for two other problems exposed 
in this volume: The Hodge conjecture, and the Birch and Swinnerton-Dyer conjecture. The conjecture 
will not succumb to only analytic nor algebraic methods. It would be necessarily a combination. 

The analytic subtlety 
consists in the influence of the Eulerian expansion away from its half plane of convergence. 
This is a matter touching the theory of resummation that involves a certain amount of 
complex analysis magic and fine analytic theory. 
The arithmetic subtlety is contained in the $\QQ$-independence of the prime frequencies $(\omega_p)$. 
By this we mean that a light perturbation of the prime frequencies in the Euler product,
\begin{equation*}
 \zeta(s)=\prod_p \left (1-e^{-\frac{\omega_p}{2\pi} s}\right )^{-1} \ ,
\end{equation*}
even assuming the subsistence of a meromorphic extension, will very likely destroy the Riemann Hypothesis (see Beurling zeta function).
Is very hard to fully exploit this simple fact in a subtle analytical argument. There is no need to say 
that there is no approach so far that  
exploits both facts in full. There hasn't been any genuine analytic-algebraic 
attack. History shows a long list of failures and partial results, hard to improve, which are genuinely 
purely analytic or purely algebraic...and all of them very far from the goal. 

Reviewing the analytic progress towards the Riemann Hypothesis is quite frustrating. The pathetic attempts to 
enlarge the ridiculous zero free region in the critical strip is a perfect example of what 
brute force can do without fully exploiting fundamental arithmetic aspects of the problem. 
And we can go on and on, with such examples, 
where the goal is becoming more and more on improving the epsilons, which evidences the fact  
that the gap towards the conjectured results will never be filled by such methods without some very 
original input.

Reviewing the algebraic progress towards the Riemann Hypothesis is equally frustrating. Only the toy 
models for 
non-transcendental functions have been deal with. 
Pushing such methods to transcendental zeta functions has been tried and
seems well out of scope. The Grothendieckian approach to divinize and solve with a trivial 
corollary, does not seem
well adapted to the analytic nature of the problem (indeed to any hard concrete analytic problem). The 
rich structure of the universe of zeta-functions, allows to combinatorially construct 
$L$-functions with a
complete lack of analytic information. This explains the proliferation of interrelated 
conjectures that after all 
is just another evidence that some analytic results are missed.

One may bet that the real hard problem is the 
Riemann Hypothesis for Riemann zeta-function, and all the vast other generalizations would fall 
from the techniques, in a pretty 
similar way, or with some extra not so hard technique.

All of this is masterly summarized in the following quotation, that every researcher in the field 
should read with humility:

\medskip

\textit {There have probably been very few attempts at proving the Riemann
hypothesis, because, simply, no one has ever had any really good idea
for how to go about it. (A. Selberg)}

\subsection{Why is it important?}

As Gauss would put it: Mathematics is the queen of Sciences, Number Theory is the queen of Mathematics, and 
the theory of prime numbers is the queen of Number Theory. To unveil the fundamental structure of central theories is what marks the progress of Science. 

But to solve the Riemann Hypothesis is central for many other reasons. It will vindicate, once more, the 
aesthetic intuition in Mathematics. Also from a practical point of view, it will prove at once all the numerous 
conditional results that exist. 

It is one of the actual principal challenges of Mathematics. These challenges, as the quadrature of the circle, the resolution of algebraic equations, Euclid's fifth postulate,...  are the only objective measure of progress in Mathematics.

Also, the resolution of the Riemann Hypothesis will very likely unlock the hidden structure in the universe 
of zeta-functions. This will open vast new fields of research.

And last, but not least, so much effort has been put on it, that the solution will save thousands of hours 
of mathematical research.

It is worth to note that the Riemann Hypothesis is the only Millenium Problem that was on Hilbert's list 
of open problem. Let's see what Hilbert and others had to say:

\medskip

\textit{...it still remains to prove the correctness of an exceedingly important statement of Riemann:
that the zero points of the function $\zeta(s)$ all have the real part $1/2$, except the well-known negative integral real zeros. (D. Hilbert, Description from the 8th Hilbert problem)}

\medskip

\textit{The failure of the Riemann Hypothesis would create havoc in the distribution of 
prime numbers. This fact alone singles out the Riemann Hypothesis as the main open 
question of prime number theory. (E. Bombieri)}

\medskip

\textit{The Riemann Hypothesis is the central problem and it implies many,
many things. One thing that makes it rather unusual in mathematics
today is that there must be over five hundred papers -somebody should go
and count- which start ``Assume the Riemann Hypothesis,...'' and the
conclusion is fantastic. And those ``conclusions'' would then become
theorems . . . With this one solution you would have proven five hundred
theorems or more at once. (P. Sarnak)}

\medskip

\subsection{Ingredients in the proof.}

From the contemplation of the panorama of zeta-functions, failed attempts and various results, one can guess 
some of the ingredients that the proof should contain. Intuition and experience 
tell us that it is very unlikely 
that two very distinct proofs may exist. There is almost a uniqueness paratheorem about 
uniqueness of proofs 
for results combining far away branches of mathematics. 

It is not very hard to guess some elements that are very likely to be present in the final proof. 

\begin{itemize}
 \item The vast universe of zeta and $L$-functions reveals an enormous hidden structure 
 that must be unveiled. This 
 analytic-algebraic structure is very likely to provide the key for the resolution of the Riemann Hypothesis.
 Also, this extra structure would explain the appearance for different $L$-functions in very different contexts:
 Simply there is  not enough place in one complex variable to admit more than one of such algebra of functions.
 
 \item  A pairing with some sort of positivity argument is very likely to be part of the key point in 
 the argument. Not only Deligne's proof for his toy model incorporates it, but various proposed 
 approaches (A. Weil, L. De Branges, physical approach...) also rely on this. The positivity may 
 only yield a weaker result, as for example the non-existence of a finite number of zeros 
 outside the critical line. For the full conjecture, more probably other ingredients will be needed.
 
 \item  A technique that allows the Euler product in the critical strip seems unavoidable. For example, a 
 proxy for the Euler product that has a global meaning. 
 
 \item  Heavy Fourier analysis mixed with complex variable must be used in a non-trivial form that incorporates 
 in full the basic arithmetic aspects. 
\end{itemize}

\subsection{Receipts for failure.}

This text grew up from lectures addressed to a broad and young mathematical audience. It may be important to give them some warnings. Most of these observations apply to most of the hard open problems that exist in Mathematics.

\begin{itemize}
 \item \textbf{Don't expect simple proofs to ever work.} It would be very naive to think otherwise.
 \item \textbf{Don't work on it unless you have very novel and powerful ideas.} Many of the best Mathematicians of all times have failed. Something more than existing techniques and tools are needed. You need a really good idea and striking new tools. Most of what you believe that is a good or novel idea is not. We hope that this text will help you to decide that. 
 \item \textbf{Don't work on it without a clear goal.} As mentioned, you must first decide if you believe the 
 conjecture or not. There is no point in trying to prove the conjecture one day and trying to disprove it the 
 next day. A clear goal is a source of strength that is needed.
 \item \textbf{Don't expect that the problem consists in resolving a single hard difficulty.} In this kind of hard problems many enemies are on your way, well hidden, and awaiting for you.
 \item \textbf{Don't work on it without studying previous attempts.} We know by now of several failed attempts, and you should learn from them in order to not repeat history again.
 \item \textbf{Don't go for it unless you have succeed in other serious problems.} ``Serious problems'' mean problems that have been open and well known for years. Before setting goals, you better check that they are realistic. Picasso didn't start his career with the Guernica. If you think that the Riemann Hypothesis will be your first major strike, you probably deserve failure.
 \item \textbf{Don't tell anyone out of your closer circle that you work on the problem.} Or you will be put aside on the freak category and will put unwanted pressure on you.
 \item \textbf{Do tell and discuss with your very best mathematical friends your work on the problem.} You will need to check very carefully the progress you make.
 \item \textbf{Don't get obsessed nor make your main goal of it.} Unless you want to ruin you mathematical career.
 \item \textbf{Don't work on it for a monetary reward.}  If you want to earn a million dollar and more, there are much simpler ways, eg. find a nice trick in finance and trading. Moreover you should know that the dollar is doomed...
\end{itemize}

Once having checked on this list...there is only one thing to say: Good luck!

\section{Appendix. Euler's $\Gamma$-function.}


For a complex number $s$ in the half plane $\Re s>0$, we define
\begin{equation*}
\Gamma(s)=\int_{0}^{+\infty} e^{-x}x^s \ \frac{dx}{x} \ .
\end{equation*} 
Then $\Gamma(1)=1$ and by integration by parts $\Gamma(s+1)=s\Gamma(s)$. In particular, $\Gamma(n)=(n-1)!$ 
for every $n \in \NN$, so the gamma function interpolates the factorial. The function $\Gamma$ can also be 
defined for all $s\in \CC$ as the limit
\begin{equation*}
\Gamma(s)=\lim_{n\rightarrow \infty} \frac{n^s n!}{s(s+1)\cdots (s+n)} \ .
\end{equation*} 
This shows that $\Gamma$ is always a non-zero function with simple poles at $0$ and 
the negative integers such that $\mathrm{Res}_{s=-n}\Gamma (s)=\frac{(-1)^n}{n!}$. 

\medskip

We have that  $1/\Gamma (s)$ is an entire function having $-\NN$ as divisor. 
Thus the order of $1/\Gamma (s)$  is $1$. We have then the Weierstrass factorization
\begin{equation*}
\frac{1}{\Gamma(s)}=s e^{\gamma s} \prod_{n=1}^{+\infty}(1+\frac{s}{n})e^{-\frac{s}{n}} \ .
\end{equation*} 
We have the complement formula
\begin{equation*}
\Gamma(s)\Gamma(1-s)=\frac{\pi}{\sin\pi z} \ ,
\end{equation*} 
observing that both sides of the equation have the same divisor, order $1$, and same equivalent for $s\to 1$.

In particular, the  evaluation at $s=1/2$ gives $ \Gamma(1/2)=\sqrt{\pi}$. 

\textit{Legendre duplication 
formula} can be obtained again by checking that both sides also have the same divisor, order $1$ and 
comparing asymptotic when $s\to 1$,
\begin{equation*}
\Gamma \left (\frac{s}{2}\right )\Gamma \left (\frac{s+1}{2} \right )=2^{1-s}\sqrt{\pi}\Gamma(s) \ .
\end{equation*}
And finally, we have \textit{Gauss formula}

\begin{equation*}
\Gamma\left (\frac{s}{n}\right )\Gamma\left (\frac{s+1}{n}\right )\cdots
\Gamma\left (\frac{s+n-1}{n}\right )=n^{\frac{1}{2}-s}(2\pi)^{\frac{n-1}{2}}\Gamma(s) \ .
\end{equation*}

\bibliographystyle{alpha}

\end{document}